\input amstex
\input amsppt.sty
\magnification=\magstep1
\hsize=32truecc
\vsize=22.2truecm
\baselineskip=16truept
\NoBlackBoxes
\TagsOnRight \pageno=1 \nologo
\def\Z{\Bbb Z}
\def\N{\Bbb N}

\def\Q{\Bbb Q}

\def\l{\left}
\def\r{\right}
\def\bg{\bigg}
\def\({\bg(}
\def\[{\bg\lfloor}
\def\){\bg)}
\def\]{\bg\rfloor}
\def\t{\text}
\def\f{\frac}

\def\p{\ (\roman{mod}\ p)}

\def\bi{\binom}
\def\eq{\equiv}

\def\ls{\leqslant}
\def\gs{\geqslant}
\def\mo{\roman{mod}}
\def\sign{\roman{sign}}

\def\ve{\varepsilon}

\def\da{\delta}

\def\Proof{\noindent{\it Proof}}

\def\Remark{\medskip\noindent{\it  Remark}}

\def\Ack{\medskip\noindent {\bf Acknowledgments}}
\hbox {Finite Fields Appl. 59 (2019), 246--283.}
\bigskip
\topmatter
\title Quadratic residues and related permutations and identities\endtitle
\author Zhi-Wei Sun\endauthor
\leftheadtext{Zhi-Wei Sun}
\rightheadtext{Quadratic residues and related permutations and identities}
\affil Department of Mathematics, Nanjing University\\
 Nanjing 210093, People's Republic of China
  \\  zwsun\@nju.edu.cn
  \\ {\tt http://math.nju.edu.cn/$\sim$zwsun}
\endaffil
\abstract Let $p$ be an odd prime. In this paper we investigate quadratic residues modulo $p$ and related permutations, congruences and identities. If $a_1<\ldots<a_{(p-1)/2}$ are all the quadratic residues modulo $p$ among $1,\ldots,p-1$, then
the list $\{1^2\}_p,\ldots,\{((p-1)/2)^2\}_p$ (with $\{k\}_p$ the least nonnegative residue of $k$ modulo $p$)
is a permutation of $a_1,\ldots,a_{(p-1)/2}$, and we show that the sign of this permutation
is $1$ or $(-1)^{(h(-p)+1)/2}$ according as $p\eq3\pmod 8$ or $p\eq7\pmod 8$, where $h(-p)$ is the class number
of the imaginary quadratic field $\Q(\sqrt{-p})$. To achieve this, we
evaluate the product $\prod_{1\ls j<k\ls(p-1)/2}(\cot\pi j^2/p-\cot\pi k^2/p)$ via Dirichlet's class number formula and Galois theory. We also obtain some new congruences and identities in product forms; for example, we determine
the exact value of $$\prod_{1\ls j<k\ls p-1}\cos\pi\frac{aj^2+bjk+ck^2}p$$ for any $a,b,c\in\Z$ with $ac(a+b+c)\not\eq0\pmod p$.
\endabstract
\thanks 2010 {\it Mathematics Subject Classification}.\,Primary 11A15, 05A05, 11T24, 33B10;
Secondary 11A07, 11R32.
\newline \indent Supported by the National Natural Science
Foundation of China (grant no. 11571162) and the NSFC-RFBR Cooperation and Exchange Program (grant no. 11811530072).
\endthanks
\endtopmatter
\document

\heading{1. Introduction}\endheading

Let $n$ be any positive integer. A permutation $\sigma$ on the set $\{1,\ldots,n\}$ is said to be odd or even according as
$$\text{Inv}(\sigma):=|\{(i,j):\ 1\ls i<j\ls n\ \text{and}\ \sigma(i)>\sigma(j)\}|$$
is odd or even. The sign of the permutation $\sigma$ is given by $\sign(\sigma)=(-1)^{\t{Inv}(\sigma)}$.
For integers $a$ and $b\not=0$ with $\gcd(b,n)=1$, we use $\{a/b\}_n$ to denote the unique integer $r\in\{0,\ldots,n-1\}$ with $a/b\eq r\pmod n$ (i.e., $a\eq br\pmod n$).

Let $p$ be an odd prime and let $a\in\Z$ with $p\nmid a$. Then $\pi_a(k)=\{ak\}_p$ with $1\ls k\ls p-1$ is a permutation on $\{1,\ldots,p-1\}$. Zolotarev's lemma (cf.
[DH] and [Z]) asserts that $\sign(\pi_a)$ coincides with
the Legendre symbol $(\f ap)$.

Frobenius (cf. [BC]) extended Zolotarev's lemma as follows: If $a\in\Z$
is relatively prime to a positive odd integer $n$, then the sign of the permutation $\pi_a(k)=\{ak\}_n$ $(0\ls k\ls n-1)$ on $\{0,\ldots,n-1\}$ equals the Jacobi symbol $(\f an)$.

Let $n>1$ be an odd integer and let $a$ be any integer relatively prime to $n$. For each $k=1,\ldots,(n-1)/2$ let $\pi^*_a(k)$ be the unique $r\in\{1,\ldots,(n-1)/2\}$ with $ak$
congruent to $r$ or $-r$ modulo $n$. For the permutation $\pi^*_a$ on $\{1,\ldots,(n-1)/2\}$, Pan [P06] showed that its sign is given by
$$\sign(\pi^*_a)=\l(\f an\r)^{(n+1)/2}.$$

Let $m>1$ be an odd integer, and let $a_1<\ldots<a_{\varphi(m)}$
be all the numbers among $1,\ldots,m-1$ relatively prime to $m$. For each $k\in\{1,\ldots,m-1\}$ with $\gcd(k,m)=1$, let $\sigma_m(k)=\bar k$ be the inverse of $k$ modulo $m$, that is,
$\bar k\in\{1,\ldots,m-1\}$ and $k\bar k\eq1\pmod m$. For $k=1,\ldots,(m-1)/2$ with $\gcd(k,m)=1$, let $\tau_m(k)$ be the unique integer $k^*\in\{1,\ldots,(m-1)/2\}$ such that $kk^*$ is congruent to $1$ or $-1$ modulo $m$. Clearly, $\sigma_m$ is a permutation of $a_1,\ldots,a_{\varphi(m)}$, and $\tau_m$ is the permutation of $a_1,\ldots,a_{\varphi(m)/2}$.
Our first theorem determines $\sign(\sigma_m)$ and $\sign(\tau_m)$.

\proclaim{Theorem 1.1} Suppose that $m=\prod_{s=1}^rp_s^{a_s}$, where $p_1,\ldots,p_r$
 are distinct odd primes and $a_1,\ldots,a_r$ are positive integers. Then we have
$$\sign(\sigma_m)=-1\iff r=1\ \t{and}\ p_1\eq1\pmod4.\tag1.1$$
Also, $\sign(\tau_m)=-1$ if and only if $r=1\ \&\ (p_1\eq1\ \t{or}\ 4a_1+3\pmod8)$, or $(r=2\ \&\ p_1+p_2\eq0\pmod4).$ In particular, when $m$ is an odd prime we have
$$\sign(\sigma_m)=-\l(\f{-1}m\r)\ \ \t{and}\ \ \sign(\tau_m)=-\l(\f 2m\r).\tag1.2$$
\endproclaim

Let $p$ be an odd prime. By Wilson's theorem,
$$(-1)^{(p-1)/2}\l(\f{p-1}2!\r)^2\eq\prod_{k=1}^{(p-1)/2}k(p-k)=(p-1)!\eq-1\pmod p.\tag1.3$$
Write $p=2n+1$ and let $a_1,\ldots,a_{n}$ be the list of all the $n$ quadratic residues among $1,\ldots,p-1$ in the ascending order. It is well known that the list
$$\{1^2\}_p,\ldots,\{n^2\}_p$$
is a permutation of $a_1,\ldots,a_{n}$.
Clearly, the sign of this permutation is just the sign of the product
$$S_p:=\prod_{1\ls i< j\ls (p-1)/2}(\{j^2\}_p-\{i^2\}_p).\tag1.4$$
(An empty product like $S_3$ is regarded to have the value $1$.)
It is easy to determine this product modulo $p$.
In fact,
$$\prod_{1\ls i<j\ls n}(j^2-i^2)\eq\cases-n!\pmod p&\t{if}\ p\eq1\pmod4,
\\1\pmod p&\t{if}\ p\eq3\pmod4,\endcases\tag1.5$$
because
$$\align \prod_{1\ls i<j\ls n}(j-i)\times\prod_{1\ls i<j\ls n}(j+i)
=&\prod_{k=1}^nk^{|\{i\gs1:\ k+i\ls n\}|}\times\prod_{k=1}^{p-1}k^{|\{1\ls i<k/2:\ k-i\ls n\}|}
\\=&\prod_{k=1}^n k^{n-k}\times\prod_{k=1}^nk^{\lfloor(k-1)/2\rfloor}(p-k)^{\lfloor k/2\rfloor}
\\\eq&(-1)^{\sum_{k=0}^n\lfloor k/2\rfloor}(n!)^{n-1}\pmod p
\endalign$$
and $(n!)^2\eq(-1)^{n+1}\pmod p$ by (1.3).
Note that if $p\eq3\pmod4$ then
$$\prod_{1\ls i<j\ls(p-1)/2}(i^2+j^2)\eq(-1)^{\lfloor(p+1)/8\rfloor}\pmod p\tag1.6$$
(cf. Problem N.2 of [Sz, pp.\,364-365]).

\medskip

Inspired by (1.5) and (1.6), we obtain the following general result.

\proclaim{Theorem 1.2} Let $p$ be an odd prime.

{\rm (i)} If $p\eq1\pmod4$, then
$$\prod\Sb 1\ls i<j\ls(p-1)/2\\p\nmid i^2+j^2\endSb (i^2+j^2)\eq(-1)^{\lfloor(p-5)/8\rfloor}\pmod p.
\tag1.7$$

{\rm (ii)} Let $a,b,c\in\Z$ with $ac(a+b+c)\not\eq0\pmod p$, and set $\Delta=b^2-4ac$. Then
$$\prod\Sb 1\ls i<j\ls p-1\\p\nmid ai^2+bij+cj^2\endSb(ai^2+bij+cj^2)\eq
\cases(\f{a(a+b+c)}p)\pmod p&\t{if}\ p\mid \Delta,\\-(\f{ac(a+b+c)\Delta}p)\pmod p&\t{if}\ p\nmid \Delta.\endcases\tag1.8$$
If $a+c=0$, then
$$\aligned&\prod^{(p-1)/2}\Sb i,j=1\\ p\nmid ai^2+bij+cj^2\endSb(ai^2+bij+cj^2)
\\\eq&
\cases\pm\f{p-1}2!\pmod p&\t{if}\ (\f{\Delta}p)=-1\ \t{or}\ (p\mid\Delta\ \&\ (\f{2b}p)=1),
\\\pm1\pmod p&\t{if}\ (\f{\Delta}p)=1\ \t{or}\ (p\mid\Delta\ \&\ (\f{2b}p)=-1).
\endcases\endaligned\tag1.9$$

{\rm (iii)} Let $a,b,c\in\Z$ with $p\nmid ac$ and $p\mid a+b+c$. Then
$$\prod\Sb 1\ls i<j\ls p-1\\p\nmid ai^2+bij+cj^2\endSb(ai^2+bij+cj^2)
\eq\cases(-1)^{N_p(a/c)}(\f{2c(a-c)}p)\pmod p&\t{if}\ p\nmid a-c,
\\(-1)^{(p+1)/2}(\f ap)\pmod p&\t{if}\ p\mid a-c,\endcases\tag1.10$$
where $N_p(x):=|\{1\ls k\ls(p-1)/2:\ \{kx\}_p>k\}|$ for any $p$-adic integer $x$.

{\rm (iv)} Let $a,b,c\in\Z$ with $p\mid ac$. Then
$$\aligned&\prod\Sb 1\ls i<j\ls p-1\\p\nmid ai^2+bij+cj^2\endSb(ai^2+bij+cj^2)
\\\eq&\cases-(\f bp)\pmod p&\t{if}\ p\mid a,\ p\nmid b\ \t{and}\ p\mid c,
\\-(\f cp)\pmod p&\t{if}\ p\mid a,\ p\nmid bc\ \t{and}\ p\mid b+c,
\\(-1)^{N_p(-c/b)}(\f 2p)\pmod p&\t{if}\ p\mid a\ \t{and}\ p\nmid bc(b+c),
\\(-1)^{(p+1)/2}(\f cp)\pmod p&\t{if}\ p\mid a,\ p\mid b\ \t{and}\ p\nmid c,
\\(-1)^{(p+1)/2}(\f ap)\pmod p&\t{if}\ p\nmid a,\ p\mid b\ \t{and}\ p\mid c,
\\(-1)^{(p+1)/2}(\f bp)\pmod p&\t{if}\ p\nmid ab,\ p\mid a+b\ \t{and}\ p\mid c,
\\(-1)^{N_p(-a/b)}(\f 2p)\pmod p&\t{if}\ p\nmid ab(a+b)\ \t{and}\ p\mid c.
\endcases\endaligned\tag1.11$$
\endproclaim

To determine the sign of $S_p$ for an arbitrary prime $p\eq3\pmod4$, we need to establish the following theorem
via Dirichlet's class number formula and Galois theory.

\proclaim{Theorem 1.3} Let $p>3$ be a prime and let $\zeta=e^{2\pi i/p}$. Let $a$ be any integer not divisible by $p$.

{\rm (i)} If $p\eq1\pmod4$, then
$$\prod_{k=1}^{(p-1)/2}(1-\zeta^{ak^2})=\sqrt p\,\ve_p^{-(\f ap)h(p)},\tag1.12$$
where $\ve_p$ and $h(p)$ are the fundamental unit and the class number of the real quadratic field
$\Q(\sqrt p)$ respectively. If $p\eq3\pmod 4$, then
$$\prod_{k=1}^{(p-1)/2}(1-\zeta^{ak^2})=(-1)^{(h(-p)+1)/2}\l(\f ap\r)\sqrt {p}\,i.\tag1.13$$

{\rm (ii)} If $p\eq1\pmod4$, then
$$\prod_{1\ls j<k\ls(p-1)/2}(\zeta^{aj^2}-\zeta^{ak^2})^2
=(-1)^{(p-1)/4}p^{(p-3)/4}\ve_p^{(\f ap)h(p)}.\tag1.14$$
When $p\eq3\pmod 4$, we have
$$\aligned&\prod_{1\ls j<k\ls(p-1)/2}(\zeta^{aj^2}-\zeta^{ak^2})
\\=&\cases(-p)^{(p-3)/8}&\t{if}\ p\eq3\pmod8,
\\(-1)^{(p+1)/8+(h(-p)-1)/2}(\f ap)p^{(p-3)/8}i&\t{if}\ p\eq7\pmod8,
\endcases\endaligned\tag1.15$$
where $h(-p)$ is the class number of the imaginary quadratic field $\Q(\sqrt{-p})$.
\endproclaim
\Remark\ 1.1. For any prime $p\eq3\pmod4$, it is known that $2\nmid h(-p)$; moreover, L. J. Mordell [M61] proved that $\f{p-1}2!\eq(-1)^{(h(-p)+1)/2}\pmod p$ if $p>3$.
In the case $a=1$, Theorem 1.3(i) appeared in [Ch]. Our proof of (1.15) utilizes the congruence (1.5).
\medskip

Since
$$\sin\pi\theta=\f i2e^{-i \pi\theta}(1-e^{2\pi i\theta})
\ \ \ \t{and}\ \ \ 2\cos\pi\theta \times \sin \pi\theta =\sin 2\pi\theta,$$
we can easily deduce the following corollary from Theorem 1.3(i).

\proclaim{Corollary 1.1} Let $p>3$ be a prime and let $a\in\Z$ with $p\nmid a$.
Then
$$\aligned&2^{(p-1)/2}\prod_{k=1}^{(p-1)/2}\sin\pi\f{ak^2}p
\\=&(-1)^{(a+1)\lfloor(p+1)/4\rfloor}\sqrt p
\times\cases\ve_p^{-(\f ap)h(p)}&\t{if}\ p\eq1\pmod4,
\\(-1)^{(h(-p)+1)/2}(\f ap)&\t{if}\ p\eq3\pmod4,\endcases
\endaligned\tag1.16$$
and
$$2^{(p-1)/2}\prod_{k=1}^{(p-1)/2}\cos\pi\f{ak^2}p=
\cases(-1)^{a(p-1)/4}\ve_p^{(1-(\f 2p))(\f ap)h(p)}&\t{if}\ p\eq1\pmod4,
\\(-1)^{(a+1)(p+1)/4}&\t{if}\ p\eq3\pmod4.\endcases\tag1.17$$
\endproclaim

For any odd prime $p$, we define
$$s(p):=\l|\l\{(j,k):\ 1\ls j<k\ls\f{p-1}2\ \t{and}\ \{j^2\}_p>\{k^2\}_p\r\}\r|\tag1.18$$
and
$$t(p):=\l|\l\{(j,k):\ 1\ls j<k\ls\f{p-1}2\ \t{and}\ \{k^2-j^2\}_p>\f p2\r\}\r|.\tag1.19$$
For example,  $s(11)=t(11)=4$  since $(\{1^2\}_{11},\ldots,\{5^2\}_{11})=(1,4,9,5,3),$
 $$\{(j,k):\ 1\ls  j < k \ls 5\ \&\ \{j^2\}_{11}> \{k^2\}_{11}\}=\{(2,5),\, (3,4),\, (3,5),\, (4,5)\},$$
 and
  $$\l\{(j,k):\ 1\ls  j < k \ls 5\ \&\ \{k^2-j^2\}_{11}> \f{11}2\r\}=\{(1,3),\, (2,5),\, (3,4),\, (4,5)\}.$$

From Theorem 1.3 we deduce the following result.

\proclaim{Theorem 1.4} Let $p$ be an odd prime.

{\rm (i)} We have
$$\sign (S_p)= (-1)^{s(p)}=(-1)^{t(p)}=\cases1&\t{if}\ p\eq3\pmod 8,
\\(-1)^{(h(-p)+1)/2}&\t{if}\ p\eq7\pmod 8.
\endcases\tag1.20$$

{\rm (ii)} Let $a\in\Z$ with $p\nmid a$. Then $$\aligned&\prod_{1\ls j<k\ls (p-1)/2}\csc\pi\f{a(k^2-j^2)}p
=\prod_{1\ls j<k\ls(p-1)/2}\l(\cot\pi\f{aj^2}p-\cot\pi\f{ak^2}p\r)
\\=&\cases (2^{p-1}/p)^{(p-3)/8}&\t{if}\ p\eq3\pmod8,
\\(-1)^{(h(-p)+1)/2}(\f ap)(2^{p-1}/p)^{(p-3)/8}&\t{if}\ p\eq7\pmod8.\endcases
\endaligned\tag1.21$$
In the case $p\eq1\pmod4$, we have
$$\aligned&(-1)^{(a-1)(p-1)/4}\prod_{1\ls j<k\ls(p-1)/2}\csc\pi\f{a(k^2-j^2)}p
\\=&\ve_p^{-(\f ap)h(p)(p-1)/2}\prod_{1\ls j<k\ls(p-1)/2}\l(\cot\pi\f{aj^2}p-\cot\pi\f{ak^2}p\r)
\\=&\pm(2^{p-1}p^{-1})^{(p-3)/8}\ve_p^{-(\f ap)h(p)/2}.
\endaligned\tag1.22$$
\endproclaim
\Remark\ 1.2.  The values of $s(p)$ for the first 2500 odd primes $p$ are available from [S18, A319311].
That $2\mid s(p)$ for any prime $p\eq3\pmod8$ might have a combinatorial proof.
\medskip

With the help of Theorem 1.4, we also get the following result.

\proclaim{Theorem 1.5} Let $p$ be an odd prime and let $\zeta=e^{2\pi i/p}$. Let $a\in\Z$ with $p\nmid a$.
Then
$$\aligned&(-1)^{a\f{p+1}2\lfloor\f{p-1}4\rfloor}2^{(p-1)(p-3)/8}\prod_{1\ls j<k\ls(p-1)/2}\cos\pi\f{a(k^2-j^2)}p
\\=&
\prod_{1\ls j<k\ls(p-1)/2}(\zeta^{aj^2}+\zeta^{ak^2})
=\cases1&\t{if}\ p\eq3\pmod4,
\\\pm\ve_p^{(\f ap)h(p)((\f2p)-1)/2}&\t{if}\ p\eq1\pmod4.\endcases
\endaligned\tag1.23$$
\endproclaim

For a real number $x$ let $\{x\}$ denote its fractional part $x-\lfloor x\rfloor$.
If $p$ is an odd prime and $1\ls j<k\ls(p-1)/2$, then
$$\align \cos2\pi\f{k^2-j^2}p<0\iff&\cos2\pi\l|\l\{\f{k^2}p\r\}-\l\{\f{j^2}p\r\}\r|<0
\\\iff& \f14<\l|\l\{\f{k^2}p\r\}-\l\{\f{j^2}p\r\}\r|<\f34.
\endalign$$
Thus Theorem 1.5 with $a=2$ yields the following corollary.

\proclaim{Corollary 1.2} For any prime $p\eq3\pmod4$, we have
$$\l|\l\{(j,k):\ 1\ls j<k\ls\f{p-1}2\ \t{and}\ \f14<\l|\l\{\f{k^2}p\r\}-\l\{\f{j^2}p\r\}\r|<\f34\r\}\r|
\eq0\pmod2.\tag1.24$$
\endproclaim

Motivated by the congruences (1.5)-(1.6) and Theorems 1.2-1.5, we establish the following theorem.

\proclaim{Theorem 1.6} Let $p$ be an odd prime.

{\rm (i)} Let $a\in\Z$ with $p\nmid a$. Then
$$\aligned&\prod\Sb 1\ls j<k\ls(p-1)/2\\p\nmid j^2+k^2\endSb\sin\pi\f{a(j^2+k^2)}p
\\=&\l(\f p{2^{p-1}}\r)^{(p-(\f{-1}p)-4)/8}\times
\cases\ve_p^{(\f ap)h(p)(1+(\f 2p))/2}&\t{if}\ p\eq1\pmod4,
\\(-1)^{(p-3)/8}&\t{if}\ p\eq3\pmod8,
\\(-1)^{(p+1)/8+(h(-p)+1)/2}(\f ap)&\t{if}\ p\eq7\pmod8.\endcases
\endaligned\tag1.25$$
Also,
$$\prod_{1\ls j<k\ls (p-1)/2}\cos\pi\f{a(j^2+k^2)}p=(-1)^{a\f{p+1}2\lfloor\f{p-1}4\rfloor}2^{-\f{p-1}2\lfloor\f{p-3}4\rfloor}
\tag1.26$$
and
$$\aligned&\prod\Sb 1\ls j<k\ls(p-1)/2\\p\nmid j^2+k^2\endSb\l(\cot\pi\f{aj^2}p+\cot\pi\f{ak^2}p\r)
\\=&(2^{p-1}p^{-1})^{(p-(\f{-1}p)-4)/8}
\times\cases\ve_p^{(\f ap)h(p)(p+(\f 2p)-4)/2}&\t{if}\ p\eq1\pmod4,
\\(-1)^{(p-3)/8}&\t{if}\ p\eq3\pmod 8,
\\(-1)^{(p+1)/8+(h(-p)+1)/2}(\f ap)&\t{if}\ p\eq7\pmod 8.\endcases
\endaligned\tag1.27$$

{\rm (ii)} Let $a,b,c\in\Z$ with $ac(a+b+c)\not\eq0\pmod p$. Set $\Delta=b^2-4ac$
and $$m=\sum_{1\ls j<k\ls p-1\atop p\mid aj^2+bjk+ck^2}(aj^2+bjk+ck^2).\tag1.28$$
Then
$$\aligned&(-1)^m\l(2^{p-1}p^{-1}\r)^{(p-3-(\f{\Delta}p))/2}\prod_{1\ls j<k\ls p-1\atop p\nmid aj^2+bjk+ck^2}\sin\pi\f{aj^2+bjk+ck^2}p
\\=&\cases(-1)^{(b+(\f{\Delta}p))\f{p-1}4}\ve_p^{h(p)((1-p+p(\f{\Delta}p)^2)(\f ap)+(\f cp)+(\f {a+b+c}p))}&\t{if}\ p\eq1\pmod4,
\\(-1)^{a+b\f{p-3}4}(\f{a(a+b+c)}p)&\t{if}\ 4\mid p-3\ \&\ p\mid\Delta,
\\(-1)^{a+(b-1)\f{p-3}4+\f{h(-p)+1}2}(\f{ac(a+b+c)\Delta}p)&\t{if}\ 4\mid p-3\ \&\ p\nmid\Delta.
\endcases
\endaligned\tag1.29$$
We also have
$$\aligned&2^{(p-1)(p-3-(\f{\Delta}p))/2}\prod_{1\ls j<k\ls p-1}\cos\pi\f{aj^2+bjk+ck^2}p
\\=&\cases(-1)^{b(p-1)/4}\ve_p^{h(p)((\f2p)-1)((1-p+p(\f{\Delta}p)^2)(\f ap)+(\f cp)+(\f{a+b+c}p))}&\t{if}\ p\eq1\pmod4,
\\(-1)^{a+b(p-3)/4+(\f{\Delta}p)(p+1)/4}&\t{if}\ p\eq3\pmod4.
\endcases\endaligned\tag1.30$$
\endproclaim
\Remark\ 1.3. Under the notation in Theorem 1.6, as $4a(aj^2+bjk+ck^2)=(2aj+bk)^2-\Delta k^2$ we have $m=0$ in the case $(\f{\Delta}p)=-1$. It seems sophisticated to determine the parity of $m$ in the case $(\f{\Delta}p)\gs0$.
\medskip

Let $p$ be any odd prime and let $a\in\Z$ with $p\nmid a$. For $1\ls j<k\ls(p-1)/2$, clearly
$$\align\{aj^2\}_p+\{ak^2\}_p>p\iff&\{aj^2\}_p>\{-ak^2\}_p
\\\iff&\cot\pi\f{aj^2}p<\cot\pi\f{-ak^2}p
\\\iff&\cot\pi\f{aj^2}p+\cot\pi\f{ak^2}p<0.
\endalign$$
Thus (1.27) yields the following consequence.

\proclaim{Corollary 1.3} Let $p$ be an odd prime and let $a\in\Z$ with $p\nmid a$. For
$$N:=\l|\l\{(j,k):\ 1\ls j<k\ls\f{p-1}2:\ \{aj^2\}_p+\{ak^2\}_p>p\r\}\r|,\tag1.31$$
we have
$$(-1)^{N}=\cases1&\t{if}\ p\eq1\pmod4,\\(-1)^{(p-3)/8}&\t{if}\ p\eq3\pmod 8,
\\(-1)^{(p+1)/8+(h(-p)+1)/2}(\f ap)&\t{if}\ p\eq7\pmod 8.
\endcases\tag1.32$$
\endproclaim

We are going to show Theorems 1.1-1.2, Theorem 1.3, Theorems 1.4-1.5 and Theorem 1.6 in Sections 2-5 respectively.
In Section 6 we pose some conjectures for further research.

\heading{2. Proofs of Theorems 1.1-1.2}\endheading

\proclaim{Lemma 2.1} Suppose that $m=\prod_{s=1}^rp_s^{a_s}$, where $p_1,\ldots,p_r$
 are distinct odd primes and $a_1,\ldots,a_r$ are positive integers. Then
$$\l|\l\{1\ls k<\f m2:\ \gcd(k,m)=1\ \t{and}\ \bar k<\f m2\r\}\r|\eq\da_{r,1}\pmod2,\tag2.1$$
where $\bar k$ is inverse of $k$ modulo $m$ $($i.e., $1\ls k\ls m-1$ and $k\bar k\eq1\pmod m)$, and $\da_{r,1}$
is $1$ or $0$ according as $r=1$ or not.
Also, the number
$$M:=\l|\l\{(i,j):\ 1\ls i<j<\f m2\ \t{and}\ ij\eq\pm1\ (\mo\ m)\r\}\r|\tag2.2$$
is odd if and only if $r=1\ \&\ (p_1\eq1\ \t{or}\ 4a_1+3\pmod8)$, or $(r=2\ \&\ p_1+p_2\eq0\pmod4).$
\endproclaim
\Proof.  By Prop. 4.2.3 of [IR, p.\,46], for each $\ve\in\{\pm1\}$ and $1\ls s\ls r$, we have
$$\align&|\{0\ls x\ls p_s^{a_s}-1:\ x^2\eq\ve\pmod{p_s^{a_s}}\}|
\\=&|\{0\ls x\ls p_s-1:\ x^2\eq\ve\pmod{p_s}\}|
\\=&\cases2&\t{if}\ \ve=1\ \t{or}\ p_s\eq1\pmod4,
\\0&\t{otherwise}.\endcases
\endalign$$
Thus, by applying the Chinese Remainder Theorem we see that
$$|\{0\ls x\ls m-1:\ x^2\eq1\pmod m\}|=2^{r}\tag2.3$$
and
$$|\{0\ls x\ls m-1:\ x^2\eq\pm1\pmod m\}|=2^{(1+\da)r},\tag2.4$$
where $\da$ is $1$ or $0$ according as whether $p_s\eq1\pmod4$ for all $s=1,\ldots,r$.

Set $n=(m-1)/2$ and
$$S=\{(k,\bar k):\ \gcd(k,m)=1\ \&\ 1\ls k,\bar k\ls n\}.$$
 Clearly, $(k,\bar k)\in S$ if and only if $(\bar k,k)\in S$. Note that
$$k=\bar k\in\{1,\ldots,n\}\iff 1\ls k\ls n\ \&\ k^2\eq1\ (\mo\ m).$$
Therefore
$$|S|\eq\l|\l\{1\ls x<\f m2:\ x^2\eq1\pmod m\r\}\r|=2^{r-1}\eq\da_{r,1}\pmod2$$
in view of (2.3). This proves (2.1).

In light of (2.4), we have
$$\align 2M=&|\{(i,j):\ 1\ls i,j\ls n\ \&\ ij\eq\pm1\ (\mo\ m)\}|
\\&-|\{1\ls x\ls n:\ x^2\eq\pm1\ (\mo\ m)\}|
\\=&|\{(i,\tau_m(i)):\ 1\ls i\ls n\ \&\ \gcd(i,m)=1\}|-2^{(1+\da)r-1}
\\=&\f{\varphi(m)}2-2^{(1+\da)r-1}=\f12\prod_{s=1}^r p_s^{a_s-1}(p_s-1)-2^{(1+\da)r-1},
\endalign$$
which implies that $M$ is odd if and only if $r=1\ \&\ (p_1\eq1\ \t{or}\ 4a_1+3\pmod8)$, or $(r=2\ \&\ p_1+p_2\eq0\pmod4).$
This concludes the proof. \qed

\medskip
\noindent{\it Proof of Theorem 1.1}. Set $n=(m-1)/2$. Clearly $\overline{m-k}=m-\bar k$ for all $1\ls k\ls m-1$ with $\gcd(k,m)=1$. If $1\ls i<j\ls m-1$ with $\gcd(i,m)=\gcd(j,m)=1$, then $m-j<m-i$ and
$$(\bar j-\bar i)(\overline{m-i}-\overline{m-j})=(\bar j-\bar i)(m-\bar i-(m-\bar j))=(\bar j-\bar i)^2>0.$$
If $1\ls i<j\ls m-1$, $\gcd(i,m)=\gcd(j,m)=1$ and $(m-j,m-i)=(i,j)$, then $1\ls i\ls n$, $j=m-i$ and $\bar j-\bar i=m-2\bar i$.
Thus
$$\sign(\sigma_m)=(-1)^{|\{1\ls i\ls n:\ \gcd(i,m)=1\ \&\ \bar i>n\}|}
=(-1)^{\varphi(m)/2-\da_{r,1}}=(-1)^{\da_{r,1}(p_1+1)/2}$$
by applying (2.1). This proves (1.1).

Now we turn to show (1.2). For $i,j\in\{1\ls k\ls n:\ \gcd(k,m)=1\}$ with $i<j$, if $i^*<j^*$ then
$$(j^*-i^*)((j^*)^*-(i^*)^*)=(j^*-i^*)(j-i)>0;$$
if $j^*<i^*$ then
$$(j^*-i^*)((i^*)^*-(j^*)^*)=(j^*-i^*)(i-j)>0;$$
if $i=i^*$ and $j=j^*$ then $j^*-i^*>0$; if $(j^*,i^*)=(i,j)$ then $j^*-i^*=i-j<0$.
In view of this, we see that
$$\sign(\tau_m)=(-1)^{|\{1\ls i\ls n:\ \gcd(i,m)=1\ \&\ i<i^*\}|}=(-1)^{M},$$
where $M$ is given by (2.2).
So the second assertion in Theorem 1.1 holds by Lemma 2.1.

The proof of Theorem 1.1 is now complete. \qed

\proclaim{Lemma 2.2} Let $p$ be an odd prime, and let $a,b,c\in\Z$ with $a$ or $b$ not divisible by $p$. Then
$$\sum_{x=0}^{p-1}\l(\f{ax^2+bx+c}p\r)=\cases-(\f ap)&\t{if}\ p\nmid b^2-4ac,
\\(p-1)(\f ap)&\t{if}\ p\mid b^2-4ac.\endcases\tag2.5$$
\endproclaim
\Remark\ 2.1. (2.5) in the case $p\mid a$ is trivial. When $p\nmid a$, (2.5) is a known result
(see, e.g., [BEW, p.\,58]).

\proclaim{Lemma 2.3} Let $p$ be any odd prime, and define
$$r(n):=\l|\l\{(j,k):\ 1\ls j<k\ls\f{p-1}2\ \t{and}\ j^2+k^2\eq n\pmod p\r\}\r|$$
for $n=0,\ldots,p-1$. Then
$$r(0)=\cases(p-1)/4&\t{if}\ p\eq1\pmod4,\\0&\t{if}\ p\eq3\pmod4.\endcases\tag2.6$$
If $n\in\{1,\ldots,p-1\}$, then
$$r(n)=\l\lfloor\f{p+1}8\r\rfloor-\f{1+(\f 2p)}2\cdot\f{1+(\f np)}2.\tag2.7$$
\endproclaim
\Proof. If $p\eq3\pmod4$, then $(\f{-1}p)=-1$ and hence $r(0)=0$. When $p\eq1\pmod4$, we have $q^2\eq-1\pmod p$
for some $q\in\Z$, and hence $r(0)=(p-1)/4$ since
$$j^2+k^2\eq0\pmod p\iff j\eq\pm qk\pmod p\iff k\eq \mp qj\pmod p.$$

Below we let $n\in\{1,\ldots,p-1\}$. Observe that
$$\align2r(n)+\f{1+(\f{2n}p)}2=&2r(n)+\l|\l\{1\ls k\ls\f{p-1}2:\ k^2+k^2\eq n\pmod p\r\}\r|
\\=&\l|\l\{(j,k):\ 1\ls j,k\ls\f{p-1}2\ \t{and}\ j^2+k^2\eq n\pmod p\r\}\r|
\\=&\l|\l\{1\ls x\ls p-1:\ \l(\f xp\r)=1\ \t{and}\ \l(\f{n-x}p\r)=1\r\}\r|
\\=&\sum_{x=1}^{p-1}\f{1+(\f xp)}2\cdot\f{1+(\f{n-x}p)}2-\f{1+(\f np)}2\cdot\f{1+(\f{n-n}p)}2
\\=&\f{p-1}4+\f14\(\sum_{x=0}^{p-1}\l(\f xp\r)+\sum_{x=0}^{p-1}\l(\f{n-x}p\r)-\l(\f np\r)\)
\\&+\f14\l(\f{-1}p\r)\sum_{x=0}^{p-1}\l(\f{x^2-nx}p\r)-\f{1+\l(\f np\r)}4
\\=&\f{p-(\f{-1}p)}4-\f{1+(\f np)}2
\endalign$$
with the help of Lemma 2.2. This yields (2.7). \qed

\proclaim{Lemma 2.4} Let $p$ be an odd prime and let $a,b,c\in\Z$ with $ac(a+b+c)\not\eq0\pmod p$.
Write $\Delta=b^2-4ac$.
For each $n=0,1\ldots,p-1$, we have
$$\aligned&|\{(j,k):\ 1\ls j<k\ls p-1\ \t{and}\ aj^2+bjk+ck^2\eq n\pmod p\}|
\\=&\cases\f12(p-3-(\f{\Delta}p)-(\f np)((1-p+p(\f{\Delta}p)^2)(\f ap)+(\f cp)+(\f{a+b+c}p)))&\t{if}\ n\not=0,
\\\f{p-1}2(1+(\f{\Delta}p))&\t{if}\ n=0.\endcases
\endaligned\tag2.8$$
\endproclaim
\Proof. Let $L$ denote the left-hand side of (2.8). Then
$$\align L=&\l|\l\{(j,k):\ 1\ls j<k\ls \f{p-1}2\ \t{and}\ aj^2+bjk+ck^2\eq n\pmod p\r\}\r|
\\&+\big|\big\{(p-j,p-k):\ 1\ls k\ls \f{p-1}2,\ k<j\ls p-1,\
\\&\qquad\quad\t{and}\ a(p-j)^2+b(p-j)(p-k)+c(p-k)^2\eq n\pmod p\big\}\big|
\\=&\big|\big\{(j,k):\ 1\ls k\ls\f{p-1}2,\ 0\ls j\ls p-1,\ j\not=0,k,\
\\&\qquad\quad\t{and}\ (2aj+bk)^2-\Delta k^2\eq4an\pmod p\big\}\big|
\endalign$$
and hence
$$\align L=&\sum_{k=1}^{(p-1)/2}\l(1+\l(\f{4an+\Delta k^2}p\r)\r)
\\&-\l|\l\{1\ls k\ls \f{p-1}2:\ (bk)^2-\Delta k^2\eq 4an\pmod p\r\}\r|
\\&-\l|\l\{1\ls k\ls \f{p-1}2:\ ((2a+b)k)^2-\Delta k^2\eq 4an\pmod p\r\}\r|.
\endalign$$
In the case $n=0$, this yields
$$L=\f{p-1}2\l(1+\l(\f{\Delta}p\r)\r).$$
When $1\ls n\ls p-1$, by the above we have
$$\align L=&\sum_{x=1}^{p-1}\f{1+(\f xp)}2\l(1+\l(\f{\Delta x+4an}p\r)\r)-\f{1+(\f{cn}p)}2-\f{1+(\f{(a+b+c)n}p)}2
\\=&\f{p-1}2+\f12\sum_{x=0}^{p-1}\l(\f xp\r)+\f12\(\sum_{x=0}^{p-1}\l(\f{\Delta x+4an}p\r)-\l(\f{4an}p\r)\)
\\&+\f12\sum_{x=0}^{p-1}\l(\f{\Delta x^2+4anx}p\r)-1-\f12\l(\f np\r)\l(\l(\f cp\r)+\l(\f{a+b+c}p\r)\r)
\\=&\f{p-3}2-\f12\l(\f{\Delta}p\r)-\f12\l(\f np\r)\l(\l(1-p\da_{(\f{\Delta}p),0}\r)\l(\f ap\r)+\l(\f cp\r)+\l(\f{a+b+c}p\r)\r)
\endalign$$
with the help of the identity
$$\sum_{x=0}^{p-1}\l(\f{\Delta x^2+4anx}p\r)=-\l(\f{\Delta}p\r)$$
from Lemma 2.2. This proves (2.8). \qed

\proclaim{Lemma 2.5} Let $p$ be an odd prime, and let $a,b,c\in\Z$ with $a+c=0$ and $abc\not\eq0\pmod p$.
Set $\Delta=b^2-4ac$. Then
$$(-1)^{|\{(i,j):\ 1\ls i,j\ls(p-1)/2\ \&\ p\mid ai^2+bij+cj^2\}|}
=\cases1&\t{if}\ (\f{\Delta}p)=-1,\\(\f2p)&\t{if}\ (\f{\Delta}p)=0,\\(\f{-1}p)&\t{if}\ (\f{\Delta}p)=1.
\endcases\tag2.9$$
\endproclaim
\Proof. Define
$$N=\l|\l\{(i,j):\ 1\ls i,j\ls\f{p-1}2\ \&\ p\mid ai^2+bij+cj^2\r\}\r|.$$

{\it Case} 1. $(\f{\Delta}p)=-1$.

In this case,
$$4a(ai^2+bij+cj^2)=(2ai+bj)^2-\Delta j^2\not\eq0\pmod p$$
for all $i,j=1,\ldots,p-1$. Thus $N=0$ and $(-1)^N=1$.

{\it Case} 2. $(\f{\Delta}p)=0$.

In this case, $p$ divides $\Delta=b^2+4a^2$, hence $(\f{-1}p)=1$ and $p\eq1\pmod 4$.
As
$$\f{b^2}{(2a)^2}\eq-1\eq\l(\f{p-1}2!\r)^2\pmod p,$$
for some $k=0,1$ we have $x:=(-1)^k\f{p-1}2!\eq -b/(2a)\pmod p$. Thus
$$\align N=&\l|\l\{(i,j):\ 1\ls i,j\ls\f{p-1}2\ \&\ i\eq jx\pmod p\r\}\r|
\\=&\f{p-1}2-\l|\l\{1\ls j\ls \f{p-1}2:\ \{jx\}_p>\f p2\r\}\r|
\endalign$$
and hence
$$(-1)^N=(-1)^{(p-1)/2}\l(\f xp\r)=\l(\f{((p-1)/2)!}p\r)=\l(\f 2p\r)$$
by using Gauss' Lemma (cf. [IR, p.\,52]) and [S19, Lemma 2.3].

{\it Case}\ 3. $(\f{\Delta}p)=1$.

In this case $\da^2\eq\Delta$ for some $\da\in\Z$ with $p\nmid\Delta$. Let $x_1$ and $x_2$
be integers with $x_1\eq(-b+\da)/(2a)\pmod p$ and $x_2\eq(-b-\da)/(2a)\pmod p$. Then
$x_1\not\eq x_2\pmod p$, $x_1+x_2\eq-b/a\pmod p$ and $x_1x_2\eq c/a=-1\pmod p$.
Thus
$$\align N=&\l|\l\{1\ls i,j\ls\f{p-1}2:\ i\eq jx_s\pmod p\ \t{for some}\ s=1,2\r\}\r|
\\=&\sum_{s=1}^2\l(\f{p-1}2-\l|\l\{1\ls j\ls \f{p-1}2:\ \{jx_s\}_p>\f p2\r\}\r|\r).
\endalign$$
Applying Gauss' Lemma we obtain that
$$(-1)^N=(-1)^{p-1}\prod_{s=1}^2\l(\f{x_s}p\r)=\l(\f{x_1x_2}p\r)=\l(\f{-1}p\r).$$

In view of the above, we have completed the proof of Lemma 2.5. \qed

\proclaim{Lemma 2.6} For any odd prime $p$, we have
$$\prod_{1\ls i<j\ls p-1}(j-i)\eq-\l(\f 2p\r)\f{p-1}2!\pmod p.\tag2.10$$
\endproclaim
\Proof. Clearly $\bi{p-1}k\eq(-1)^k\pmod p$ for all $k=0,\ldots,p-1$. Also, $(p-1)!\eq-1\pmod p$ by Wilson's theorem. Thus
$$\align\prod_{1\ls i<j\ls p-1}(j-i)=&\prod_{j=2}^{p-1}(j-1)!=\prod_{k=1}^{p-2}k!
\\=&\f{p-1}2!\prod_{0<k<(p-1)/2}\f{(p-1)!}{\bi{p-1}k}
\eq\f{p-1}2!\prod_{0<k<(p-1)/2}(-1)^{k-1}
\\\eq&\f{p-1}2!(-1)^{(p-3)(p-5)/8}=\f{p-1}2!(-1)^{(p^2-9)/8}
\\\eq&-\l(\f 2p\r)\f{p-1}2!\pmod p.
\endalign$$
This proves (2.10). \qed

\proclaim{Lemma 2.7} Let $p$ be an odd prime, and let $a,b\in\Z$ with $a\not\eq0,1\pmod p$. Then
$$|\{x\in\{0,1\ldots,p-1\}:\ \{ax+b\}_p>x\}|=\f{p-1}2.\tag2.11$$
\endproclaim
\Proof. For $x\in\{0,\ldots,p-1\}$, obviously
$$\{ax+b+1\}_p>x\iff p-1>\{ax+b\}_p\gs x.$$
As $a\not\eq0,1\pmod p$, we have
$$|\{x\in\{0,\ldots,p-1\}:\ \{ax+b\}_p=x\}|=1=|\{x\in\{0,\ldots,p-1\}:\ \{ax+b\}_p=p-1\}|.$$
Thus
$$|\{x\in\{0,\ldots,p-1\}:\ \{ax+b+1\}_p>x\}|=|\{x\in\{0,\ldots,p-1\}:\ \{ax+b\}_p>x\}|.$$

In view of the above, it suffices to prove (2.11) for $b=0$.
For $x=1,\ldots,p-1$, clearly
$\{ax\}_p\not=x$ (since $a\not\eq1\pmod p$), and also
$$\{ax\}_p>x\iff p-\{ax\}_p<p-x\iff \{a(p-x)\}_p<p-x.$$
So (2.11) holds for $b=0$. This concludes the proof. \qed

\medskip
\noindent{\it Proof of Theorem 1.2}. (i) Let $r(n)$ be as in Lemma 2.3. Then
$$\align \prod\Sb 1\ls i<j\ls(p-1)/2\\p\nmid i^2+j^2\endSb(i^2+j^2)
\eq&\prod_{n=1}^{p-1}n^{r(n)}=((p-1)!)^{\lfloor(p+1)/8\rfloor}\prod^{p-1}\Sb n=1\\(\f np)=1\endSb
n^{-(1+(\f2p))/2}
\\\eq&(-1)^{\lfloor(p+1)/8\rfloor}\prod_{k=1}^{(p-1)/2}(k^2)^{-(1+(\f 2p))/2}
\\\eq&(-1)^{\lfloor(p+1)/8\rfloor}\l((-1)^{(p+1)/2}\r)^{(1+(\f 2p))/2}\pmod p
\endalign$$
with the help of (1.3). This yields (1.7) if $p\eq1\pmod4$. It also proves (1.6) in the case $p\eq3\pmod4$.

(ii) If $p\mid \Delta$, then by Lemma 2.4 and (1.3) we have
$$\align&\prod\Sb 1\ls i<j\ls(p-1)/2\\ p\nmid ai^2+bij+cj^2\endSb(ai^2+bij+cj^2)
\\\eq&\prod_{n=1}^{p-1}n^{\f{p-3}2+\f{1-p}2(\f ap)+\f12((\f cp)+(\f{a+b+c}p))-\f12(1+(\f np))((1-p)(\f ap)+(\f cp)+(\f{a+b+c}p))}
\\\eq&(-1)^{\f{p-3}2+\f{1-p}2(\f ap)+\f12((\f cp)+(\f{a+b+c}p))}\prod_{k=1}^{(p-1)/2}(k^2)^{-((1-p)(\f ap)+(\f cp)+(\f{a+b+c}p))}
\\\eq&(-1)^{\f12((\f cp)+(\f{a+b+c}p))-1}\l((-1)^{(p+1)/2}\r)^{(1-p)(\f ap)+(\f cp)+(\f{a+b+c}p)}
\\=&(-1)^{\f12((\f cp)+(\f{a+b+c}p))-1}=\l(\f{c(a+b+c)}p\r)=\l(\f{a(a+b+c)}p\r)\pmod p
\endalign$$
since $4ac\eq b^2\pmod p$.
Similarly, when $p\nmid\Delta$, by Lemma 2.4 and (1.3) we have
$$\align&\prod\Sb 1\ls i<j\ls(p-1)/2\\ p\nmid ai^2+bij+cj^2\endSb(ai^2+bij+cj^2)
\\\eq&\prod_{n=1}^{p-1}n^{\f{p-3}2+\f12((\f ap)+(\f cp)+(\f{a+b+c}p)-(\f{\Delta}p))-\f12(1+(\f np))((\f ap)+(\f cp)+(\f{a+b+c}p))}
\\\eq&(-1)^{\f{p-3}2+\f12((\f ap)+(\f cp)+(\f{a+b+c}p)-(\f{\Delta}p))}\prod_{k=1}^{(p-1)/2}(k^2)^{-((\f ap)+(\f cp)+(\f{a+b+c}p))}
\\\eq&(-1)^{\f{p-3}2+\f12((\f ap)+(\f cp)+(\f{a+b+c}p)-(\f{\Delta}p))}\l((-1)^{(p+1)/2}\r)^{(\f ap)+(\f cp)+(\f{a+b+c}p)}\pmod p\endalign$$
and hence
$$\align\prod\Sb 1\ls i<j\ls(p-1)/2\\ p\nmid ai^2+bij+cj^2\endSb(ai^2+bij+cj^2)
\eq&(-1)^{\f12((\f ap)+(\f cp))}(-1)^{\f12((\f{a+b+c}p)-(\f{\Delta}p))}
\\=&-\l(\f{ac}p\r)\l(\f{(a+b+c)\Delta}p\r)\pmod p.
\endalign$$
This proves (1.8).

Now assume that $a+c=0$. We deduce (1.9) from (1.8).
Observe that
$$\align &\prod\Sb 1\ls i<j\ls p-1\\p\nmid ai^2+bij+cj^2\endSb(ai^2+bij+cj^2)
\\=&\prod^{(p-1)/2}\Sb i,j=1\\p\nmid ai^2-bij+cj^2\endSb(ai^2+bi(p-j)+c(p-j)^2)\times
\prod\Sb 1\ls i<j\ls (p-1)/2\\p\nmid ai^2+bij+cj^2\endSb(ai^2+bij+cj^2)
\\&\times\prod\Sb1\ls j<i\ls (p-1)/2\\p\nmid ai^2+bij+cj^2\endSb(a(p-i)^2+b(p-i)(p-j)+c(p-j)^2)
\\\eq&\prod^{(p-1)/2}\Sb i,j=1\\p\nmid ai^2-bij+cj^2\endSb(ai^2-bij+cj^2)\times\prod^{(p-1)/2}\Sb i,j=1\\p\nmid ai^2+bij+cj^2\endSb(ai^2+bij+cj^2)
\bigg/\prod_{i=1}^{(p-1)/2}(ai^2+bi^2+ci^2)
\\\eq&\f{(-1)^m}{((p-1)/2)!^2}\l(\f {a+b+c}p\r)\prod^{(p-1)/2}\Sb i,j=1\\p\nmid ci^2+bij+aj^2\endSb(ci^2+bij+aj^2)\times
\prod^{(p-1)/2}\Sb i,j=1\\p\nmid ai^2+bij+cj^2\endSb(ai^2+bij+cj^2)
\\=&\f{(-1)^{m}}{((p-1)/2)!)^2}\l(\f {a+b+c}p\r)
\prod^{(p-1)/2}\Sb i,j=1\\p\nmid ai^2+bij+cj^2\endSb(ai^2+bij+cj^2)^2\pmod p,
\endalign$$
where
$$\align m=&\l|\l\{(i,j):\ 1\ls i,j\ls\f{p-1}2\ \t{and}\ \ p\nmid ai^2-bij+cj^2\r\}\r|
\\=&\l|\l\{(i,j):\ 1\ls i,j\ls\f{p-1}2\ \t{and}\ \ p\nmid ci^2+bij+aj^2\r\}\r|.\endalign$$
Combining this with (1.8) and noting $ac=-a^2$, we see that
$$\prod^{(p-1)/2}\Sb i,j=1\\p\nmid ai^2+bij+cj^2\endSb(ai^2+bij+cj^2)^2
\eq(-1)^m\l(\f{p-1}2!\r)^2\times\cases(\f ap)\pmod p&\t{if}\ p\mid\Delta,\\-(\f{-\Delta}p)\pmod p&\t{if}\ p\nmid\Delta.
\endcases\tag2.12$$

By Lemma 2.5,
$$(-1)^{((p-1)/2)^2-m}=\cases1&\t{if}\ (\f{\Delta}p)=-1,\\(\f 2p)&\t{if}\ (\f{\Delta}p)=0,
\\(\f{-1}p)&\t{if}\ (\f{\Delta}p)=1.\endcases$$
If $p$ divides $\Delta=b^2+4a^2$, then $(\f{-1}p)=1$, $(-1)^m=(\f 2p)$, $(\f b{2a})^2\eq-1\eq(\f{p-1}2!)^2\pmod p$ and
$$\l(\f{b}p\r)=\l(\f{\pm2\times((p-1)/2)!}p\r)\l(\f ap\r)=\l(\f ap\r)$$
with the help of [S19, Lemma 2.3].
Thus, in view of (2.12) and (1.3), we have (1.9) if $p\mid\Delta$.
When $(\f{\Delta}p)=-1$, we have
$$-(-1)^m\l(\f{-\Delta}p\r)=(-1)^m\l(\f{-1}p\r)=1$$
and hence (1.9) holds in view of (2.12). If $(\f{\Delta}p)=1$, then
$$-(-1)^m\l(\f{-\Delta}p\r)=-(-1)^m\l(\f{-1}p\r)=-\l(\f{-1}p\r)\eq\f1{((p-1)/2)!^2}\pmod p$$
and hence (1.9) follows from (2.12).

(iii) If $a\eq c\pmod p$, then $b\eq-a-c\eq-2a\pmod p$ and hence
$$\align&\prod_{1\ls i<j\ls p-1\atop p\nmid ai^2+bij+cj^2}(ai^2+bij+cj^2)
\\\eq&\prod_{1\ls i<j\ls p-1}a(j-i)^2=a^{\f{p-1}2(p-2)}\prod_{1\ls i<j\ls p-1}(j-i)^2
\\\eq&\l(\f ap\r)\l(\f{p-1}2!\r)^2\eq(-1)^{(p+1)/2}\l(\f ap\r)\pmod p
\endalign$$
with the help of (2.10) and (1.3).

Now assume that $a\not \eq c\pmod p$. For $1\ls i<j\ls p-1$, clearly
$$ai^2+bij+cj^2\eq (i-j)(ai-cj)\eq c(j-i)\l(j-\f aci\r)\pmod p.$$
Note that
$$\align\prod_{1\ls i<j\ls p-1\atop p\nmid ai-cj}\l(j-\f aci\r)
\eq&\prod_{r=1}^{p-1}r^{|\{(i,j):\ 1\ls i<j\ls p-1\ \&\ j-\f aci\eq r\pmod p\}|}
\\\eq&\prod_{r=1}^{p-1}r^{|\{1\ls i\ls p-1:\ \{r+\f aci\}_p>i\}|}
\\\eq&(p-1)!^{(p-1)/2-1}\eq(-1)^{(p+1)/2}\pmod p
\endalign$$
with the help of Lemma 2.7 and Wilson's theorem. Also,
$$\align\prod_{1\ls i<j\ls p-1\atop p\mid ai-cj}c(j-i)\eq&\prod_{i=1\atop\{\f aci\}_p>i}^{p-1}c\l(\f aci-i\r)
\\\eq&\prod_{i=1\atop\{\f aci\}_p>i}^{(p-1)/2}(a-c)i\times\prod_{i=1\atop\{\f ac(p-i)\}_p>p-i}^{(p-1)/2}(a-c)(p-i)
\\\eq&\prod_{i=1}^{(p-1)/2}(a-c)i\times (-1)^{|\{1\ls i\ls\f{p-1}2:\ \{\f aci\}_p<i\}|}
\\\eq&\l(\f {a-c}p\r)\f{p-1}2!(-1)^{(p-1)/2-N_p(a/c)}\pmod p
\endalign$$
and hence
$$\align \prod_{1\ls i<j\ls p-1\atop p\nmid ai-cj}c(j-i)
\eq&\f{\prod_{1\ls i<j\ls p-1}c(j-i)}{(\f {a-c}p)\f{p-1}2!(-1)^{(p-1)/2-N_p(a/c)}}
\\\eq&\f{-c^{(p-1)(p-2)/2}(\f 2p)\f{p-1}2!}{(\f {a-c}p)\f{p-1}2!(-1)^{(p-1)/2-N_p(a/c)}}
\\\eq&\l(\f{2c(a-c)}p\r)(-1)^{(p+1)/2+N_p(a/c)}\pmod p
\endalign$$
with the help of (2.10).
Therefore
$$\align\prod_{1\ls i<j\ls p-1\atop p\nmid ai^2+bij+cj^2}(ai^2+bij+cj^2)
\eq&\prod_{1\ls i<j\ls p-1\atop p\nmid ai-cj}c(j-i)\times \prod_{1\ls i<j\ls p-1\atop p\nmid ai-cj}\l(j-\f aci\r)
\\\eq&\l(\f{2c(a-c)}p\r)(-1)^{N_p(a/c)}\pmod p.
\endalign$$
This proves part (iii) of Theorem 1.2.

(iv) In the spirit of our proof of Theorem 1.2(iii), we may show Theorem 1.2(iv).
To illustrate this, here we handle the case $p\mid a$ and $p\nmid bc(b+c)$ in details.
Note that
$$\prod_{1\ls i<j\ls p-1\atop p\nmid ai^2+bij+cj^2}(ai^2+bij+cj^2)
\eq\prod_{1\ls i<j\ls p-1\atop p\nmid bi+cj}cj\times \prod_{1\ls i<j\ls p-1\atop p\nmid bi+cj}\l(j+\f bci\r)\pmod p.\tag2.13$$
Similar to the second paragraph in (iii), we have
$$\prod_{1\ls i<j\ls p-1\atop p\nmid bi+cj}\l(j+\f bci\r)\eq(-1)^{(p+1)/2}\pmod p.\tag2.14$$
Observe that
$$\align\prod_{1\ls i<j\ls p-1\atop p\mid bi+cj}cj
=&\prod_{j=1\atop\{-\f cbj\}_p<j}^{(p-1)/2}cj\times\prod_{j=1\atop\{-\f cb(p-j)\}_p<p-j}^{(p-1)/2}c(p-j)
\\\eq&\prod_{j=1}^{(p-1)/2}cj\times (-1)^{|\{1\ls j\ls\f{p-1}2:\ \{-\f cbj\}_p>j\}|}
\\\eq&\l(\f cp\r)\f{p-1}2!(-1)^{N_p(-c/b)}\pmod p
\endalign$$
and hence
$$\align\prod_{1\ls i<j\ls p-1\atop p\nmid bi+cj}cj\eq&\f{\prod_{j=1}^{(p-1)/2}(cj)^{j-1}(c(p-j))^{p-j-1}}
{(\f cp)\f{p-1}2!(-1)^{N_p(-c/b)}}
\\\eq&\prod_{j=1}^{(p-1)/2}\f{(-1)^j(cj)^{p-1}}{cj}\times\f{(\f cp)(-1)^{N_p(-c/b)}}{\f{p-1}2!}
\\\eq&(-1)^{(p^2-1)/8+N_p(-c/b)}\l(\f{p-1}2!\r)^{-2}
\\\eq&(-1)^{(p+1)/2+N_p(-c/b)}\l(\f 2p\r)
\pmod p
\endalign$$
with the help of (1.3). Combining this with (2.13) and (2.14), we obtain that
$$\prod_{1\ls i<j\ls p-1\atop p\nmid ai^2+bij+cj^2}(ai^2+bij+cj^2)
\eq (-1)^{N_p(-c/b)}\l(\f 2p\r)\pmod p.$$
This ends our proof. \qed

\heading{3. Proof of Theorem 1.3}\endheading

To prove Theorem 1.3, we need some known results.

\proclaim{Lemma 3.1} Let $p$ be an odd prime, and let $\zeta=e^{2\pi i/p}$.

{\rm (i)} For any $a\in\Z$ with $p\nmid a$, we have
$$\prod_{n=1}^{p-1}(1-\zeta^{an})=p\tag3.1$$
and
$$\sum_{x=0}^{p-1}\zeta^{ax^2}=\l(\f ap\r)\sqrt{(-1)^{(p-1)/2}p}.\tag3.2$$

{\rm (ii)} If $p\eq1\pmod4$, then
$$\prod_{n=1}^{p-1}(1-\zeta^n)^{(\f np)}=\ve_p^{-2h(p)}.\tag3.3$$

{\rm (iii)} When $p\eq3\pmod4$, we have
$$ph(-p)=-\sum_{k=1}^{p-1}k\l(\f kp\r),\tag3.4$$
and also
$$\l|\l\{1\ls k\ls \f{p-1}2:\ \l(\f kp\r)=-1\r\}\r|\eq\f{h(-p)+1}2\pmod2\tag3.5$$
provided $p>3$.
\endproclaim
\Remark\ 3.1. This lemma is well known. For any $a\in\Z$ with $p\nmid a$, we have (3.1) since
$$\prod_{n=1}^{p-1}(1-\zeta^{an})=\prod_{k=1}^{p-1}(1-\zeta^k)=\lim_{x\to1}\f{x^{p}-1}{x-1}=p,$$
and we have (3.2) by Gauss' evaluation of quadratic Gauss sums (cf. [IR, pp.\,70-75]).
Part (ii) and the first assertion in part (iii) are Dirichlet's class number formula.
The second assertion in part (iii) was pointed out by Mordell [M61].
\medskip

\proclaim{Lemma 3.2} Let $p$ be an odd prime and let $n\in\{1,\ldots,p-1\}$.
Then
$$\aligned&\l|\l\{(j,k):\ 1\ls j,k\ls\f{p-1}2\ \t{and}\ j^2-k^2\eq n\pmod p\r\}\r|
\\=&\l\lfloor\f{p-1}4\r\rfloor-\cases1&\t{if}\ p\eq1\pmod4\ \t{and}\ (\f np)=1,
\\0&\t{otherwise}.\endcases\endaligned\tag3.6$$
\endproclaim
\Proof. Let $L$ denote the left-hand side of (3.6). Then
$$\align
L=&\l|\l\{1\ls x\ls p-1:\ \l(\f xp\r)=1\ \t{and}\ \l(\f{n+x}p\r)=1\r\}\r|
\\=&\sum_{x=1}^{p-1}\f{(\f xp)+1}2\cdot\f{(\f{x+n}p)+1}2-\f{(\f{p-n}p)+1}2\cdot\f{(\f {p-n+n}p)+1}2
\\=&\f{p-1}4+\f14\sum_{x=0}^{p-1}\l(\f xp\r)+\f14\(\sum_{x=0}^{p-1}\l(\f{x+n}p\r)-\l(\f np\r)\)
\\&+\f14\sum_{x=0}^{p-1}\l(\f{x(x+n)}p\r)-\f{(\f{-n}p)+1}4
\\=&\f{p-3-(\f np)-(\f{-n}p)}4
\endalign$$
with the help of Lemma 2.2.
This yields (3.6). \qed

\medskip
\noindent{\it Proof of Theorem 1.3}. Let $\varphi_a$ be the element of the Galois group
$\t{Gal}(\Q(\zeta)/\Q)$ with $\varphi_a(\zeta)=\zeta^a$.
In view of (3.2),
$$\varphi_a\l(\sqrt{(-1)^{(p-1)/2}p}\r)=\varphi_a\(\sum_{x=0}^{p-1}\zeta^{x^2}\)=\sum_{x=0}^{p-1}\zeta^{ax^2}=\l(\f ap\r)\sqrt{(-1)^{(p-1)/2}p}.\tag3.7$$

 (i) We first handle the case $p\eq1\pmod4$.
Combining (3.1) and (3.3), we get
$$\prod^{p-1}\Sb n=1\\(\f np)=1\endSb(1-\zeta^n)^{2}
=\prod_{n=1}^{p-1}(1-\zeta^n)^{1+(\f np)}
=p\ve_p^{-2h(p)}.$$
Note that
$$\prod^{p-1}\Sb n=1\\(\f np)=1\endSb(1-\zeta^n)
=\prod^{(p-1)/2}\Sb n=1\\(\f np)=1\endSb(1-\zeta^n)(1-\zeta^{p-n})
=\prod^{(p-1)/2}\Sb n=1\\(\f np)=1\endSb|1-\zeta^n|^2>0.$$
Therefore
$$\prod_{k=1}^{(p-1)/2}(1-\zeta^{k^2})=\prod^{p-1}\Sb n=1\\(\f np)=1\endSb(1-\zeta^n)=\sqrt{p}\,\ve_p^{-h(p)}.\tag3.8$$
This proves (1.12) for $a=1$.

Write $\ve_p=u_p+v_p\sqrt{p}$ with $u_p,v_p\in\Q$. In view of (3.7),
$$\varphi_a(\ve_p)=u_p+\l(\f ap\r)v_p\sqrt{p}=\f{N(\ve_p)}{u_p-(\f ap)v_p\sqrt p}
=\cases\ve_p&\t{if}\ (\f ap)=1,\\N(\ve_p)\ve_p^{-1}&\t{if}\ (\f ap)=-1,\endcases$$
where $N(\ve_p)$ is the norm of $\ve_p$ with respect to the field extension $\Q(\sqrt p)/\Q$.
Thus, by using (3.7) and (3.8) we obtain
$$\align\prod_{k=1}^{(p-1)/2}(1-\zeta^{ak^2})
=&\varphi_a\(\prod_{k=1}^{(p-1)/2}(1-\zeta^{k^2})\)=\varphi_a\l(\sqrt p\,\ve_p^{-h(p)}\r)
\\=&\cases \sqrt{p}\,\ve_p^{-h(p)}&\t{if}\ (\f ap)=1,
\\-\sqrt p\,N(\ve_p)^{-h(p)}\ve_p^{h(p)}&\t{if}\ (\f ap)=-1.\endcases
\endalign$$
This proves (1.12) since $N(\ve_p)^{h(p)}=-1$ (cf. [Co, p.\,185 and p.\,187]).

Now we consider the case $p\eq3\pmod4$. In view of (3.4),
$$\align ph(-p)=&-\sum_{r=1}^{(p-1)/2}\l(r\l(\f rp\r)+(p-r)\l(\f{p-r}p\r)\r)
\\=&-2\sum_{r=1}^{(p-1)/2}r\l(\f rp\r)+p\sum_{r=1}^{(p-1)/2}\l(\f rp\r)
\endalign$$
and hence $p\mid\sum_{r=1}^{(p-1)/2}r(\f rp)$. Let $N=|\{1\ls r\ls(p-1)/2:\ (\f rp)=-1\}|$.
Observe that
$$\align \prod_{k=1}^{(p-1)/2}(1-\zeta^{k^2})=&\prod_{k=1}^{(p-1)/2}(1-\zeta^{(2k)^2})
\\=&\prod^{(p-1)/2}\Sb r=1\\(\f rp)=1\endSb(1-\zeta^{4r})\times\prod^{(p-1)/2}\Sb r=1\\(\f rp)=-1\endSb(1-\zeta^{4(p-r)})
\\=&\prod^{(p-1)/2}\Sb r=1\\(\f rp)=1\endSb(1-\zeta^{4r})\times\prod^{(p-1)/2}\Sb r=1\\(\f rp)=-1\endSb\f{\zeta^{4r}-1}{\zeta^{4r}}
\\=&(-1)^N\zeta^{-4\sum_{0<r<p/2,\,(\f rp)=-1}r}\prod_{r=1}^{(p-1)/2}\zeta^{2r}(\zeta^{-2r}-\zeta^{2r})
\\=&(-1)^N\zeta^{\sum_{r=1}^{(p-1)/2}2r(\f rp)}\prod_{r=1}^{(p-1)/2}(\zeta^{-2r}-\zeta^{2r})
\\=&(-1)^N\prod_{r=1}^{(p-1)/2}(\zeta^{-2r}-\zeta^{2r})
\endalign$$
and hence
$$(-1)^N\prod_{k=1}^{(p-1)/2}(1-\zeta^{k^2})=(-1)^{(p-1)/2}\prod_{r=1}^{(p-1)/2}(\zeta^{2r}-\zeta^{-2r}).\tag3.9$$

By Prop. 6.4.3 of [IR, p.\,74],
$$\prod_{r=1}^{(p-1)/2}(\zeta^{2r-1}-\zeta^{-(2r-1)})=\sqrt p\,i.$$
Thus
$$\align\sqrt p\,i\prod_{r=1}^{(p-1)/2}(\zeta^{2r}-\zeta^{-2r})=&\prod_{k=1}^{p-1}(\zeta^k-\zeta^{-k})
=\zeta^{\sum_{k=1}^{p-1}k}\prod_{k=1}^{p-1}(1-\zeta^{-2k})=p
\endalign$$
with the help of (3.1).
Combining this with (3.9) we obtain
$$(-1)^N\prod_{k=1}^{(p-1)/2}(1-\zeta^{k^2})=\sqrt p\,i.$$
Therefore, by using (3.7) we get
$$(-1)^N\prod_{k=1}^{(p-1)/2}(1-\zeta^{ak^2})=\varphi_a(\sqrt p\,i)=\l(\f ap\r)\sqrt p\, i.$$
This yields (1.13) since
$N\eq(h(-p)+1)/2\pmod2$ by (3.5).

(ii)  Observe that
$$\align &\prod_{1\ls j<k\ls (p-1)/2}(\zeta^{aj^2}-\zeta^{ak^2})^2
\\=&(-1)^{\bi{(p-1)/2}2}\prod_{1\ls j<k\ls(p-1)/2}(\zeta^{aj^2}-\zeta^{ak^2})(\zeta^{ak^2}-\zeta^{aj^2})
\\=&(-1)^{\bi{(p-1)/2}2}\prod_{k=1}^{(p-1)/2}\prod^{(p-1)/2}\Sb j=1\\j\not=k\endSb(\zeta^{ak^2}-\zeta^{aj^2})
\\=&(-1)^{\bi{(p-1)/2}2}\prod_{k=1}^{(p-1)/2}(\zeta^{ak^2})^{(p-3)/2}
\times\prod_{k=1}^{(p-1)/2}\,\prod^{(p-1)/2}\Sb j=1\\j\not=k\endSb(1-\zeta^{a(j^2-k^2)})
\\=&(-1)^{(p-1)(p-3)/8}\zeta^{\f{p-3}2\sum_{k=1}^{(p-1)/2}ak^2}\prod^{(p-1)/2}\Sb j,k=1\\j\not=k\endSb
(1-\zeta^{a(j^2-k^2)}).
\endalign$$
Clearly,
$$\sum_{k=1}^{(p-1)/2}k^2=\f{p^2-1}{24}p\eq0\pmod p.\tag3.10$$
So, with the help of Lemma 3.2, from the above we obtain
$$\align\prod_{1\ls j<k\ls (p-1)/2}(\zeta^{aj^2}-\zeta^{ak^2})^2
=&(-1)^{(p-1)(p-3)/8}\prod_{n=1}^{p-1}(1-\zeta^{an})^{\lfloor(p-1)/4\rfloor}
\\&\times\cases\prod_{0<n<p,\ (\f np)=1}(1-\zeta^{an})^{-1}&\t{if}\ p\eq1\pmod4,\\1&\t{if}\ p\eq3\pmod4.
\endcases\endalign$$
Noting (3.1) and Theorem 1.3(i) we get
$$\aligned&\prod_{1\ls j<k\ls (p-1)/2}(\zeta^{aj^2}-\zeta^{ak^2})^2
\\=&(-1)^{(p-1)(p-3)/8}p^{(p-3)/4}\times\cases\ve_p^{(\f ap)h(p)}&\t{if}\ p\eq1\pmod4,\\1&\t{if}\ p\eq3\pmod4.\endcases\endaligned\tag3.11$$
Thus (1.14) holds when $p\eq1\pmod4$.

Below we suppose $p\eq3\pmod4$ and want to show (1.15).
By (3.11) with $a=1$, for some $\ve\in\{\pm1\}$, we have
$$\prod_{1\ls j<k\ls (p-1)/2}(\zeta^{j^2}-\zeta^{k^2})=\ve (\sqrt p\,i)^{(p-3)/4}.$$
In view of (3.7), this yields that
$$\prod_{1\ls j<k\ls (p-1)/2}(\zeta^{aj^2}-\zeta^{ak^2})
=\ve\varphi_a(\sqrt p\,i)^{(p-3)/4}
=\ve \l(\l(\f ap\r)\sqrt p\,i\r)^{(p-3)/4}.\tag3.12$$
As $i^{(p-3)/4}=(-1)^{(p-7)/8}i$ if $p\eq7\pmod8$, we obtain (1.15) from (3.12) provided that
$$\ve=\cases1&\t{if}\ p\eq3\pmod8,
\\(-1)^{(h(-p)+1)/2}&\t{if}\ p\eq7\pmod8.\endcases\tag3.13$$

Now it remains to show (3.13). By (3.12), for any $r=1,\ldots,(p-1)/2$ we have
$$\prod_{1\ls j<k\ls(p-1)/2}\l(\zeta^{r^2j^2}-\zeta^{r^2k^2}\r)=\ve(\sqrt p\,i)^{(p-3)/4};$$
on the other hand,
$$\prod_{1\ls j<k\ls(p-1)/2}\l(\zeta^{r^2j^2}-\zeta^{r^2k^2}\r)
=\zeta^{r^2\sum_{1\ls j<k\ls(p-1)/2}j^2}\prod_{1\ls j<k\ls(p-1)/2}(1-\zeta^{r^2(k^2-j^2)}).$$
Combining these and noting (3.10) and (1.13), we find that
$$\align\l(\ve(\sqrt p\,i)^{(p-3)/4}\r)^{(p-1)/2}=&\prod_{1\ls j<k\ls(p-1)/2}\prod_{r=1}^{(p-1)/2}(1-\zeta^{(k^2-j^2)r^2})
\\=&\prod_{1\ls j<k\ls(p-1)/2}\l((-1)^{(h(-p)+1)/2}\l(\f{k^2-j^2}p\r)\sqrt p\,i\r).
\endalign$$
Therefore
$$\ve^{(p-1)/2}=(-1)^{\f{h(-p)+1}2\cdot\f{(p-1)(p-3)}8}\prod_{1\ls j<k\ls(p-1)/2}\l(\f{k^2-j^2}p\r)
=(-1)^{\f{h(-p)+1}2\cdot\f{p-3}4}$$
with the help of (1.5). This proves the desired (3.13) since $\ve=\ve^{(p-1)/2}$.

 The proof of Theorem 1.3 is now complete. \qed

 \heading{4. Proofs of Theorems 1.4 and 1.5}\endheading

 \proclaim{Lemma 4.1} Let $p$ be any odd prime. Then
 $$\sum_{1\ls j<k\ls(p-1)/2}(j^2+k^2)\eq \cases p\pmod {2p}&\t{if}\ p\eq5\pmod 8,
 \\0\pmod{2p}&\t{otherwise}.\endcases\tag4.1$$
 \endproclaim
 \Proof. Since
 $$\align&2\sum_{1\ls j<k\ls(p-1)/2}(j^2+k^2)+\sum_{k=1}^{(p-1)/2}(k^2+k^2)
 \\=&\sum_{j=1}^{(p-1)/2}\sum_{k=1}^{(p-1)/2}(j^2+k^2)=(p-1)\sum_{k=1}^{(p-1)/2}k^2,
 \endalign$$
 we have
 $$\sum_{1\ls j<k\ls(p-1)/2}(j^2+k^2)=\f{p-3}2\sum_{k=1}^{(p-1)/2}k^2=\f{p-3}2\cdot\f{p^2-1}{24}p\eq0\pmod p.$$
 Note that
 $$\f{p-3}2\cdot\f{p^2-1}{24}\eq\cases(p-1)/4\pmod2&\t{if}\ p\eq1\pmod4,
 \\0\pmod2&\t{if}\ p\eq3\pmod4.\endcases$$
 Therefore (4.1) holds. \qed

\medskip
\noindent{\it Proof of Theorem 1.4}.
For $1\ls j<k\ls (p-1)/2$, clearly
$$\{j^2\}_p>\{k^2\}_p\iff\cot \pi\f{j^2}p-\cot\pi\f{k^2}p<0$$
and $$\{k^2-j^2\}_p>\f p2\iff\sin2\pi\f{k^2-j^2}p<0\iff \csc2\pi\f{k^2-j^2}p<0.$$
So (1.20) follows from (1.21) and we only need to show part (ii) of Theorem 1.4.

As (1.21) holds trivially for $p=3$, below we assume $p>3$.

For $1\ls j<k\ls(p-1)/2$, clearly
 $$\align\sin\pi\f{a(k^2-j^2)}p=&\f{e^{i\pi a(k^2-j^2)/p}-e^{-i\pi a(k^2-j^2)/p}}{2i}
 \\=&\f i2e^{-i\pi a(k^2+j^2)/p}(e^{2\pi iaj^2/p}-e^{2\pi iak^2/p}).
 \endalign$$
 Combining this with Lemma 4.1, we see that
 $$\aligned&\prod_{1\ls j<k\ls (p-1)/2}\sin\pi\f{a(k^2-j^2)}p
 \\=&(-1)^{a\f{p+1}2\lfloor\f{p-1}4\rfloor}\l(\f i2\r)^{(p-1)(p-3)/8}\prod_{1\ls j<k\ls(p-1)/2}(e^{2\pi iaj^2/p}-e^{2\pi iak^2/p}).
 \endaligned\tag4.2$$

 For any real numbers $\theta_1,\theta_2\not\in\Z$, clearly
 $$\cot\pi\theta_1-\cot\pi\theta_2=\f{\cos\pi\theta_1}{\sin\pi\theta_1}-\f{\cos\pi\theta_2}{\sin\pi\theta_2}
 =\f{\sin\pi(\theta_2-\theta_1)}{\sin\pi\theta_1\sin\pi\theta_2}.$$
 Thus
 $$\align&\prod_{1\ls j<k\ls(p-1)/2}\f{\sin\pi a(k^2-j^2)/p}{\cot\pi aj^2/p-\cot\pi ak^2/p}
 \\=&\prod_{1\ls j<k\ls(p-1)/2}\sin\pi\f{aj^2}p\sin\pi\f{ak^2}p
 =\prod_{k=1}^{(p-1)/2}\l(\sin\pi\f{ak^2}p\r)^{|\{1\ls j\ls(p-1)/2:\ j\not=k\}|}
 \endalign$$
 and hence by (1.16) we have
 $$\aligned&\prod_{1\ls j<k\ls(p-1)/2}\sin\pi\f{a(k^2-j^2)}p\bigg/\prod_{1\ls j<k\ls(p-1)/2}\l(\cot\pi\f{aj^2}p-\cot\pi\f{ak^2}p\r)
 \\=&\l(\f p{2^{p-1}}\r)^{(p-3)/4}\times\cases(-1)^{(a-1)(p-1)/4}\ve_p^{-(\f ap)\f{p-3}2h(p)}&\t{if}\ p\eq1\pmod4,
 \\1&\t{if}\ p\eq3\pmod4.\endcases\endaligned\tag4.3$$
 So it suffices to determine $\prod_{1\ls j<k\ls (p-1)/2}\sin\pi a(k^2-j^2)/p$.

 {\it Case} 1. $p\eq3\pmod4$.

 In this case, by combining (4.2) and (1.15) we get
 $$\align&\prod_{1\ls j<k\ls(p-1)/2}\sin\pi\f{a(k^2-j^2)}p
 \\=&\l(\f p{2^{p-1}}\r)^{(p-3)/8}\times\cases1&\t{if}\ p\eq3\pmod 8,
 \\(-1)^{(h(-p)+1)/2}(\f ap)&\t{if}\ p\eq7\pmod 8.\endcases
 \endalign$$
 Thus (1.21) is valid with the help of (4.3).

 {\it Case} 2. $p\eq1\pmod4$.

 In this case, combining (4.2) with (1.14) we obtain
 $$\prod_{1\ls j<k\ls(p-1)/2}\sin^2\pi\f{a(k^2-j^2)}p=\l(\f p{2^{p-1}}\r)^{(p-3)/4}\ve_p^{(\f ap)h(p)}\tag4.4$$
 and hence
 $$\prod_{1\ls j<k\ls(p-1)/2}\csc\pi\f{a(k^2-j^2)}p=\pm(2^{p-1}p^{-1})^{(p-3)/8}\ve_p^{-(\f ap)h(p)/2}.$$
 In view of (4.3), we have
 $$\align &\prod_{1\ls j<k\ls(p-1)/2}\l(\sin\pi\f{a(k^2-j^2)}p\r)\l(\cot\pi\f{aj^2}p-\cot\f{ak^2}p\r)
 \\=&(-1)^{(a-1)(p-1)/4}(2^{p-1}p^{-1})^{(p-3)/4}\ve_p^{(\f ap)\f{p-3}2h(p)}
 \prod_{1\ls j<k\ls(p-1)/2}\sin^2\pi\f{a(k^2-j^2)}p.
 \endalign$$
 Combining this with (4.4) we immediately get the first equality in (1.22).

The proof of Theorem 1.4 is now complete. \qed
\medskip

\noindent{\it Proof of Theorem 1.5}. (1.23) is trivial for $p=3$. Below we assume $p>3$. In view of (4.2),
$$\align \prod_{1\ls j<k\ls(p-1)/2}\l(2\cos\pi\f{a(k^2-j^2)}p\r)
=&\prod_{1\ls j<k\ls(p-1)/2}\f{\sin \pi (2a)(k^2-j^2)/p}{\sin \pi a(k^2-j^2)/p}
\\=&(-1)^{a\f{p+1}2\lfloor\f{p-1}4\rfloor}\prod_{1\ls j<k\ls(p-1)/2}(\zeta^{aj^2}+\zeta^{ak^2}).
\endalign$$
So we have the first equality in (1.23). On the other hand, by Theorem 1.4(ii) we have
$$\prod_{1\ls j<k\ls(p-1)/2}\f{\csc\pi a(k^2-j^2)/p}{\csc\pi(2a)(k^2-j^2)/p}
=\cases1&\t{if}\ p\eq3\pmod4,\\\pm\ve_p^{(\f ap)h(p)((\f 2p)-1)/2}&\t{if}\ p\eq1\pmod4.\endcases$$
Therefore (1.23) holds.

The proof of Theorem 1.5 is now complete. \qed

\heading{5. Proof of Theorem 1.6}\endheading

\proclaim{Lemma 5.1} Let $p$ be an odd prime. Then
$$\f1p\sum\Sb 1\ls j<k\ls(p-1)/2\\p\mid j^2+k^2\endSb(j^2+k^2)\eq\cases1\pmod2&\t{if}\ p\eq5\pmod8,
\\0\pmod2&\t{otherwise}.\endcases\tag5.1$$
\endproclaim
\Proof. If $p\eq3\pmod 4$, then $(\f{-1}p)=-1$ and $j^2+k^2\not\eq0\pmod p$ for any $j,k=1,\ldots,(p-1)/2$.
So (5.1) is trivial in the case $p\eq3\pmod4$.

Now assume that $p\eq1\pmod4$. Then $q^2\eq-1\pmod p$ for some $q\in\Z$. For each $j=1,\ldots,(p-1)/2$
let $j_*$ be the unique integer $r\in\{1,\ldots,(p-1)/2\}$ with $qj$ congruent to $r$ or $-r$ modulo $p$.
Clearly, $\{j_*:\ 1\ls j\ls(p-1)/2\}=\{1,\ldots,(p-1)/2\}$. Thus
$$\align\sum\Sb 1\ls j<k\ls(p-1)/2\\p\mid j^2+k^2\endSb(j^2+k^2)
=&\f12\sum^{(p-1)/2}\Sb j,k=1\\p\mid j^2+k^2\endSb(j^2+k^2)
\\=&\f12\sum_{j=1}^{(p-1)/2}(j^2+j_*^2)=\sum_{k=1}^{(p-1)/2}k^2=\f{p^2-1}{24}p
\endalign$$
and hence
$$\f1p\sum\Sb 1\ls j<k\ls(p-1)/2\\p\mid j^2+k^2\endSb(j^2+k^2)=\f{p^2-1}{24}
\eq\f{p^2-1}8\eq\f{p-1}4\pmod2.$$
Therefore (5.1) holds. \qed

\medskip
\noindent{\it Proof of Theorem 1.6(i)}.
Let $\zeta=e^{2\pi i/p}$. As
$$\sum\Sb 1\ls j<k\ls(p-1)/2\\p\nmid j^2+k^2\endSb(j^2+k^2)\eq0\pmod{2p}$$
by Lemmas 4.1 and 5.1, we have
$$\align\prod\Sb 1\ls j<k\ls(p-1)/2\\p\nmid j^2+k^2\endSb\sin\pi\f{a(j^2+k^2)}p
=&\prod\Sb 1\ls j<k\ls(p-1)/2\\p\nmid j^2+k^2\endSb\f{-e^{-i\pi a(j^2+k^2)/p}}{2i}(1-\zeta^{a(j^2+k^2)})
\\=&\l(\f i2\r)^{|\{(j,k):\ 1\ls j<k\ls (p-1)/2\ \&\ p\nmid j^2+k^2\}|}f(a),
\endalign$$
where $$f(a):=\prod_{n=1}^{p-1}(1-\zeta^{an})^{r(n)}$$
with $r(n)$ defined as in Lemma 2.3. Note that
$$\aligned&\l|\l\{(j,k):\ 1\ls j<k\ls \f{p-1}2\ \&\ p\nmid j^2+k^2\r\}\r|
\\=&\bi{(p-1)/2}2-r(0)=\f{p-1}2\l\lfloor\f{p-3}4\r\rfloor
\endaligned\tag5.2$$
with the help of (2.6). So
$$\prod\Sb 1\ls j<k\ls(p-1)/2\\p\nmid j^2+k^2\endSb\sin\pi\f{a(j^2+k^2)}p=\l(\f i2\r)^{\f{p-1}2\lfloor\f{p-3}4\rfloor}f(a).\tag5.3$$
By (2.7), (3.1) and Theorem 1.3(i), we have
$$\align f(a)=&p^{\lfloor(p+1)/8\rfloor}\prod^{p-1}\Sb n=1\\(\f np)=1\endSb(1-\zeta^{an})^{-(1+(\f 2p))/2}
\\=&\cases p^{\lfloor(p+1)/8\rfloor}&\t{if}\ p\eq3,5\pmod 8,
\\p^{\lfloor(p+1)/8\rfloor-1/2}\ve_p^{(\f ap)h(p)}&\t{if}\ p\eq1\pmod 8,
\\(-1)^{(h(-p)-1)/2}(\f ap)p^{\lfloor(p+1)/8\rfloor-1/2}i&\t{if}\ p\eq7\pmod 8.
\endcases\endalign$$
Combining this with (5.3) we immediately get (1.25).

In light of (1.25) and (5.2),
$$\align\prod\Sb 1\ls j<k\ls(p-1)/2\\p\nmid j^2+k^2\endSb\cos\pi\f{a(j^2+k^2)}p
=&\prod\Sb 1\ls j<k\ls(p-1)/2\\p\nmid j^2+k^2\endSb\f{\sin\pi(2a)(j^2+k^2)/p}{2\sin\pi a(j^2+k^2)/p}
\\=&2^{-|\{(j,k):\ 1\ls j<k\ls(p-1)/2\ \&\ p\nmid j^2+k^2\}|}
\\=&2^{-\f{p-1}2\l\lfloor\f{p-3}4\r\rfloor}.
\endalign$$
Note also that
$$\align&\prod\Sb 1\ls j<k\ls(p-1)/2\\p\mid j^2+k^2\endSb\cos\pi\f{a(j^2+k^2)}p
\\=&(-1)^{\sum_{1\ls j<k\ls(p-1)/2\ \&\ p\mid j^2+k^2}a(j^2+k^2)/p}
=(-1)^{a\f{p+1}2\lfloor\f{p-1}4\rfloor}
\endalign$$ by Lemma 5.1.
Therefore (1.26) holds.

Observe that
$$\prod\Sb 1\ls j<k\ls(p-1)/2\\p\nmid j^2+k^2\endSb\l(\cot\pi\f{aj^2}p+\cot\pi\f{ak^2}p\r)
=\prod\Sb 1\ls j<k\ls(p-1)/2\\p\nmid j^2+k^2\endSb\f{\sin\pi a(j^2+k^2)/p}{(\sin\pi aj^2/p)(\sin\pi ak^2/p)}$$
and
$$\prod\Sb 1\ls j<k\ls(p-1)/2\\p\nmid j^2+k^2\endSb\l(\sin\pi\f{aj^2}p\r)\l(\sin\pi\f{ak^2}p\r)
=\prod_{k=1}^{(p-1)/2}\l(\sin\pi\f{ak^2}p\r)^{(p-(\f{-1}p)-4)/2}.$$
Combining these with (1.25) and (1.16), we obtain the desired (1.27).
This concludes the proof of Theorem 1.6(i). \qed

\proclaim{Lemma 5.2} Let $p>3$ be a prime and let $a,b,c\in\Z$ with $p\nmid a$. Then
$$\sum_{1\ls j<k\ls p-1}(aj^2+bjk+ck^2)\eq0\pmod p\tag5.4$$
and also $$\f1p\sum_{1\ls j<k\ls p-1}(aj^2+bjk+ck^2)\eq a\f{p-1}2+b\f{(p-1)(p-3)}8\pmod2.\tag5.5$$
\endproclaim
\Proof. Let $\Delta=b^2-4ac$. In view of (3.10), we have
$$\align &\sum_{1\ls j<k\ls p-1}(aj^2+bjk+ck^2)
\\=&\sum_{1\ls j<k\ls(p-1)/2}(aj^2+bjk+ck^2)
\\&+\sum_{1\ls k\ls(p-1)/2}\sum_{k<j\ls p-1}(a(p-j)^2+b(p-j)(p-k)+c(p-k)^2)
\\\eq&\sum_{k=1}^{(p-1)/2}\(\sum_{j=0}^{p-1}(aj^2+bjk+ck^2)-ck^2-(a+b+c)k^2\)
\\\eq&\sum_{k=1}^{(p-1)/2}\sum_{j=0}^{p-1}\f1{4a}\l((2aj+bk)^2-\Delta k^2\r)
\eq\sum_{k=1}^{(p-1)/2}\f1{4a}\sum_{r=0}^{p-1}r^2\eq0\pmod p.
\endalign$$
This proves (5.4).

Observe that
$$\align &\sum_{1\ls j<k\ls p-1}(aj^2+bjk+ck^2)
\\\eq&\sum_{1\ls j<k\ls p-1}(aj+ck)+b\sum_{1\ls j<k\ls(p-1)/2}(2j-1)(2k-1)
\\\eq&\sum_{k=1}^{p-1}\(\sum_{0<j<k} aj+ck(k-1)\)+b\bi{(p-1)/2}2
\\\eq&\sum_{k=1}^{p-1}\f a2(k^2-k)+b\f{(p-1)(p-3)}8\pmod2
\endalign$$
and hence
$$\align &\sum_{1\ls j<k\ls p-1}(aj^2+bjk+ck^2)-b\f{(p-1)(p-3)}8
\\\eq&\f a2\l(\f{(p-1)p(2p-1)}6-\f{(p-1)p}2\r)=a\f{p(p-1)(p-2)}6\eq a\f{p-1}2\pmod2.
\endalign$$
Therefore (5.5) also holds. \qed

\medskip
\noindent{\it Proof of Theorem 1.6(ii)}. By Lemma 2.4,
$$\aligned&\l|\l\{(j,k):\ 1\ls j<k\ls p-1\ \t{and}\ p\nmid aj^2+bjk+ck^2\r\}\r|
\\=&\bi{p-1}2-\f{p-1}2\l(1+\l(\f{\Delta}p\r)\r)=\f{p-1}2\l(p-3-\l(\f{\Delta}p\r)\r).
\endaligned\tag5.6$$
Let $\zeta=e^{2\pi i/p}$.
Then
$$\align&\prod_{1\ls j<k\ls p-1\atop p\nmid aj^2+bjk+ck^2}\sin\pi\f{aj^2+bjk+ck^2}p
\\=&\prod_{1\ls j<k\ls p-1\atop p\nmid aj^2+bjk+ck^2}\f{-e^{-i\pi(aj^2+bjk+ck^2)/p}}{2i}(1-\zeta^{aj^2+bjk+ck^2})
\\=&\l(\f i2\r)^{\f{p-1}2(p-3-(\f{\Delta}p))}(-1)^{a(p-1)/2+b(p-1)(p-3)/8-m}
\\&\times\prod_{1\ls j<k\ls p-1\atop p\nmid aj^2+bjk+ck^2}(1-\zeta^{aj^2+bjk+ck^2})
\endalign$$
with the help of Lemma 5.2.
In view of Lemma 2.4, (3.1) and Theorem 1.3(i), we have
$$\align &\prod_{1\ls j<k\ls p-1\atop p\nmid aj^2+bjk+ck^2}(1-\zeta^{aj^2+bjk+ck^2})
\\=&\f{\prod_{n=1}^{p-1}(1-\zeta^n)^{(p-3-(\f{\Delta}p)+(1-p+p(\f{\Delta}p)^2)(\f ap)+(\f cp)+(\f{a+b+c}p))/2}}
{\prod_{n=1}^{p-1}(1-\zeta^n)^{((1-p+p(\f{\Delta}p)^2)(\f ap)+(\f cp)+(\f{a+b+c}p))(1+(\f np))/2}}
\\=&\f{p^{(p-3-(\f{\Delta}p)+(1-p+p(\f{\Delta}p)^2)(\f ap)+(\f cp)+(\f{a+b+c}p))/2}}
{\prod_{k=1}^{(p-1)/2}(1-\zeta^{k^2})^{(1-p+p(\f{\Delta}p)^2)(\f ap)+(\f cp)+(\f{a+b+c}p)}}
\\=&p^{(p-3-(\f{\Delta}p)/2}
\\&\times
\cases \ve_p^{h(p)((1-p+p(\f{\Delta}p)^2)(\f ap)+(\f cp)+(\f{a+b+c}p))}&\t{if}\ p\eq1\pmod4,
\\((-1)^{(h(-p)-1)/2}i)^{(1-p+p(\f{\Delta}p)^2)(\f ap)+(\f cp)+(\f{a+b+c}p)}&\t{if}\ p\eq3\pmod4.
\endcases
\endalign$$
Therefore
$$\align &(-1)^m\prod_{1\ls j<k\ls p-1\atop p\nmid aj^2+bjk+ck^2}\sin\pi\f{aj^2+bjk+ck^2}p
\\=&(-1)^{a(p-1)/2+b(p-1)(p-3)/8}i^{\f{p-1}2(p-3-(\f{\Delta}p))}\l(\f p{2^{p-1}}\r)^{(p-3-(\f{\Delta}p))/2}
\\&\times
\cases \ve_p^{h(p)((1-p+p(\f{\Delta}p)^2)(\f ap)+(\f cp)+(\f{a+b+c}p))}&\t{if}\ p\eq1\pmod4,
\\(-1)^{(\f{\Delta}p)\f{h(-p)-1}2+\f{1-p}2(\f ap)+\f12((\f cp)+(\f{a+b+c}p))}
i^{p(\f{\Delta}p)^2(\f ap)}&\t{if}\ p\eq3\pmod4.
\endcases
\endalign$$
It is easy to see that this implies (1.29).

Clearly,
$$\prod_{1\ls j<k\ls p-1\atop p\mid aj^2+bjk+ck^2}\cos\pi\f{aj^2+bjk+ck^2}p
=\prod_{1\ls j<k\ls p-1\atop p\mid aj^2+bjk+ck^2}(-1)^{(aj^2+bjk+ck^2)/p}=(-1)^m.$$
On the other hand, by (5.6) we have
$$\align&2^{\f{p-1}2(p-3-(\f{\Delta}p))}\prod_{1\ls j<k\ls p-1\atop p\nmid aj^2+bjk+ck^2}\cos\pi\f{aj^2+bjk+ck^2}p
\\=&\prod_{1\ls j<k\ls p-1\atop p\nmid aj^2+bjk+ck^2}\f{\sin\pi(2aj^2+2bjk+2ck^2)/p}{\sin\pi(aj^2+bjk+ck^2)/p}.
\endalign$$
Combining these with (1.29) we immediately obtain the desired (1.30).

In view of the above, we have completed the proof of Theorem 1.6(ii). \qed

\heading{6. Some conjectures}\endheading

We are unable to determine the parities of $s(p)$ and $t(p)$ (defined by (1.18) and (1.19))
for a general prime $p\eq1\pmod4$. However, in contrast with (1.20), we formulate the following conjecture.

\proclaim{Conjecture 6.1} For any prime $p\eq1\pmod4$, we have
$$s(p)+t(p)\eq\l|\l\{1\ls k<\f p4:\ \l(\f kp\r)=1\r\}\r|\pmod 2.\tag6.1$$
\endproclaim

For any positive odd number $n$ and integer $k$, we let $R(k,n)$ denote the unique $r\in\{0,\ldots,(n-1)/2\}$
with $k$ congruent to $r$ or $-r$ modulo $n$. For example,
$$R(1^2,11)=1,\ R(2^2,11)=4,\ R(3^2,11)=2,\ R(4^2,11)=5,\ R(5^2,11)=3.$$

Motivated by Theorem 1.4 and Corollary 1.3, we pose the following conjecture.

\proclaim{Conjecture 6.2} Let $p$ be an odd prime, and let $a\in\Z$ with $p\nmid a$. Then
$$\l|\l\{(i,j): \ 1\ls i<j\ls\f{p-1}2\ \t{and}\ R(ai^2,p)>R(aj^2,p)\r\}\r|\eq\l\lfloor\f{p+1}8\r\rfloor\pmod2,\tag6.2$$
and
$$\aligned&(-1)^{|\{(i,j):\ 1\ls i<j\ls(p-1)/2\ \&\ R(ai^2,p)+R(aj^2,p)>p/2\}|}
\\=&\cases(-1)^{|\{1\ls k<\f p4:\ (\f kp)=-1\}|}(\f ap)^{(1-(\f 2p))/2}&\t{if}\ p\eq1\pmod4,
\\1&\t{if}\ p\eq3\pmod4.
\endcases\endaligned\tag6.3$$
\endproclaim
\Remark\ 6.1. We have verified (6.2) and (6.3) with $a=1$ for all odd primes $p<20000$.
\medskip

\proclaim{Conjecture 6.3} Let $p>3$ be a prime.
If $p\eq3\pmod4$, then
$$\aligned&(-1)^{|\{(j,k):\ 1\ls j<k\ls(p-1)/2\ \&\ \{j(j+1)/2\}_p>\{k(k+1)/2\}_p\}|}
\\&\qquad=(-1)^{\f{h(-p)+1}2+|\{1\ls k\ls\lfloor\f{p+1}8\rfloor:\ (\f kp)=1\}|}.\endaligned\tag6.4$$
Also,
$$\aligned&(-1)^{|\{(j,k):\ 1\ls j<k\ls(p-1)/2\ \&\ \{j(j+1)/2\}_p+\{k(k+1)/2\}_p>p\}|}
 \\=&\cases(-1)^{(p-1)/8}&\t{if}\ p\eq1\pmod8,
 \\(-1)^{|\{1\ls k<\f p4:\ (\f kp)=-1\}|}&\t{if}\ p\eq5\pmod 8,
 \\(-1)^{\f{h(-p)+1}2+|\{1\ls k\ls\lfloor\f{p+1}8\rfloor:\ (\f kp)=-1\}|}&\t{if}\ p\eq3\pmod4.
 \endcases\endaligned\tag6.5$$
\endproclaim

\proclaim{Conjecture 6.4} Let $p$ be an odd prime. If $p\eq3\pmod4$, then
$$(-1)^{\l|\l\{(j,k): \ 1\ls j<k\ls(p-1)/2\ \t{and}\ \{j(j+1)\}_p>\{k(k+1)\}_p\r\}\r|}=(-1)^{\lfloor(p+1)/8\rfloor}.\tag6.6$$
Also,
$$\aligned&(-1)^{|\{(j,k):\ 1\ls j<k\ls(p-1)/2\ \&\ \{j(j+1)\}_p+\{k(k+1)\}_p>p\}|}
\\=&\cases(-1)^{\lfloor(p-1)/8\rfloor}&\t{if}\ p\eq1\pmod4,
\\(-1)^{(h(-p)+1)/2}&\t{if}\ p>3\ \&\ p\eq3\pmod8,
\\1&\t{if}\ p\eq7\pmod8.\endcases\endaligned\tag6.7$$
\endproclaim

\proclaim{Conjecture 6.5} {\rm (i)} For any prime $p\eq5\pmod6$, we have
$$\l|\l\{1\ls k\ls\f{p-1}2:\ \{k^3\}_p>\f p2\r\}\r|-\f{p+1}6\in \{2n:\ n=0,1,2,\ldots\}\tag6.8$$
and
$$|\{(j,k):\ 1\ls j<k\ls p-1\ \t{and}\ \{j^3\}_p>\{k^3\}_p\}|\eq\f{p+1}6\pmod2.\tag6.9$$

{\rm (ii)} For any integer $m>1$, we have
$$\l|\l\{1\ls k\ls\f{p-1}2:\ \{k^m\}_p>\f p2\r\}\r|\sim \f p4\tag6.10$$
as $p\to\infty$, where $p$ is an odd prime.
\endproclaim
\Remark\ 6.2. Let $p$ be a prime with $p\eq5\pmod6$. The list $\{1^3\}_p,\ldots,\{(p-1)^3\}_p$
is a permutation of $1,\ldots,p-1$, for, if $1\ls j<k\ls p-1$ then
$$j^3-k^3=(j-k)(j^2+jk+k^2)=\f{j-k}4((2j+k)^2+3k^2)\not\eq0\pmod p.$$
See [S18, A320044] for some data related to (6.8).
Note that (6.8) implies (6.9) since for any $1\ls j<k\ls p-1$ we have $1\ls p-k<p-j\ls p-1$ and
$$(\{j^3\}_p-\{k^3\}_p)(\{(p-k)^3\}_p-\{(p-j)^3\}_p)>0.$$

\proclaim{Conjecture 6.6} Let $p$ be an odd prime. Then
$$\aligned&\l|\l\{(j,k):\ 1\ls j<k\ls \f{p-1}2\ \t{and}\ \{j^4\}_p>\{k^4\}_p\r\}\r|
\\\eq&\l\lfloor\f{p+1}8\r\rfloor+\cases(h(-p)+1)/2\pmod2&\t{if}\ p\eq7\pmod8,
\\0\pmod2&\t{otherwise}.\endcases
\endaligned\tag6.11$$
Also,
$$\aligned&\l|\l\{(j,k):\ 1\ls j<k\ls \f{p-1}2\ \t{and}\ \{j^8\}_p>\{k^8\}_p\r\}\r|
\\\eq&\cases|\{1\ls k<\f p4:\ (\f kp)=1\}|\pmod2&\t{if}\ p\eq1\pmod8,
\\0\pmod2&\t{if}\ p\eq3\pmod 8,\\(p-5)/8\pmod2&\t{if}\ p\eq5\pmod8,
\\(h(-p)+1)/2\pmod2&\t{if}\ p\eq7\pmod8.\endcases
\endaligned\tag6.12$$
\endproclaim
\Remark\ 6.3. See [S18, A309012, A319882, A319894 and A319903] for related data or similar conjectures.
\medskip

The following conjecture is motivated by Theorem 1.5.

\proclaim{Conjecture 6.7} Let $p$ be a prime with $p\eq1\pmod4$, and let $\zeta=e^{2\pi i/p}$.
Let $a$ be an integer not divisible by $p$. Then
$$\aligned&(-1)^{|\{1\ls k<p/4:\ (\f kp)=-1\}|}\prod_{1\ls j<k\ls(p-1)/2}(\zeta^{aj^2}+\zeta^{ak^2})
\\=&\cases1&\t{if}\ p\eq1\pmod8,
\\(\f ap)\ve_p^{-(\f ap)h(p)}&\t{if}\ p\eq5\pmod8.\endcases
\endaligned\tag6.13$$
\endproclaim
\Remark\ 6.4.  By K. S. Williams and J. D. Currie [WC],  for any prime $p\eq1\pmod8$ we have
$$2^{(p-1)/4}\eq(-1)^{|\{1\ls k<p/4:\ (\f kp)=-1\}|}\pmod p.$$
\medskip

The author [S19] studied the determinants of the matrices
$$\l[\l(\f{i^2+j^2}p\r)\r]_{1\ls i,j\ls(p-1)/2}
\ \ \t{and}\ \ \l[\l(\f{i^2+j^2}p\r)\r]_{0\ls i,j\ls(p-1)/2},$$
where $p$ is an odd prime. Now we conclude this section with a conjecture involving determinants.

\proclaim{Conjecture 6.8} Let $n>1$ be an odd integer. Then
$$\det[R(i^2j^2,n)]_{1\ls i,j\ls(n-1)/2}\not=0\tag6.14$$
if and only if $n$ is a prime congruent to $3$ modulo $4$.
Also,
$$\det\l[\l\lfloor \f{i^2j^2}n\r\rfloor\r]_{1\ls i,j\ls(n-1)/2}\not=0\tag6.15$$
if and only if $n$ is either $9$ or a prime greater than $7$ and congruent to $3$ modulo $4$.
\endproclaim

\Ack. The author would like to thank the two referees for their helpful comments.

\enddocument